\documentclass[11pt,letterpaper]{article}
\usepackage{amsmath, amsfonts, amsthm, amssymb, verbatim}
\usepackage{graphicx}
\usepackage[mathcal]{euscript}
\usepackage{amstext}
\usepackage{amssymb,mathrsfs}
\usepackage{bbm}
\usepackage{amsmath, amsfonts, amsthm}
\usepackage{amssymb,mathrsfs, verbatim}
\usepackage{graphicx,amstext}
\usepackage{amsmath}
\usepackage{mathtools}
\usepackage{amssymb}
\usepackage{latexsym}
\usepackage{ mathrsfs }
\usepackage{enumerate}
\usepackage{xcolor}
\usepackage[utf8]{inputenc}
\usepackage{aeguill}
\newcommand{\R}{\mathbb R}

\newcommand{\ve}{\varepsilon}
\newcommand{\vp}{\varphi}

\newcommand{\dive}{{\rm div}}

\newcommand{\p}{\partial}

\newcommand{\curl}{\,{\rm curl}}

\newcommand{\diam}{\,{\rm diam}}

\newtheorem{theorem}{Theorem}[section]
\newtheorem{proposition}[theorem]{Proposition}
 \newtheorem{remark}[theorem]{Remark}
\newtheorem{lemma}[theorem]{Lemma}
\newtheorem{corollary}[theorem]{Corollary}
\newtheorem{definition}[theorem]{Definition}

\begin{document}
\title{An analysis of the 2-D isentropic Euler Equations
for a generalized polytropic gas law}

\author{Talita Mello\footnote{Instituto de Matem\'atica, Universidade Federal do Rio de Janeiro, 
Brazil. 
E-mail:  {\sl  souzamello.talita@gmail.com}},
Wladimir Neves\footnote{Instituto de Matem\'atica, Universidade Federal do Rio de Janeiro, 
Brazil. E-mail:  {\sl wladimir@im.ufrj.br}
}
}

\date{}

\maketitle

\noindent \textit{ {\bf Key words and phrases:} 
Two-dimensional isentropic Euler equations, system of equations of mixed type, generalized polytropic gas law, quasi-potential flows.}

\vspace{0.3cm}

 \noindent {\bf MSC2010 subject classification:} 35L65, 35M30, 35Q31, 35Q35, 76N15.  

\begin{abstract} 
In this paper we develop an analysis of the compressible, isentropic Euler equations 
in two spatial dimensions for a generalized polytropic gas law. The main focus is rotational flows
in the subsonic regimes, described through the framework of 
the Euler equations expressed in self-similar variables and pseudo-velocities.
A Bernoulli type equation is derived, serving as a cornerstone for establishing a self-similar 
system tailored to rotational flows. 
We also develop an Ellipticity Principle for generalized polytropic gases, 
which is applied twice in this paper. 
To the best of the authors’ knowledge, both applications appear for the first time.
In particular, the analysis of the potential flow problem in the 
pseudo-subsonic regime is nontrivial for generalized polytropic 
gases when $\gamma < 1$. In this setting, refined techniques, 
such as the Moser iteration method combined with 
suitable a priori estimates, are required.
In the final section, the study extends to an analysis of 
a perturbed model, introducing the concept of quasi-potential flows, 
offering insights into their behavior and implications.

\end{abstract}

\tableofcontents

%%%%%%%%%%%%%%%%%%
\section{Introduction}
\label{Intro}
%%%%%%%%%%%%%%%%%%

In this paper we develop an analysis of the compressible, isentropic Euler equations 
in two spatial dimensions for a generalized polytropic gas law without assuming irrotational 
(also called potential) flow. 
Let us recall that, the compressible Euler equations is the core example 
of systems of conservation laws in multiple spatial variables, which continues to be, to this day, uncharted territory. 
However, the particular case of two spatial dimensions (2-D for short), besides being physically relevant, has recently proven to be 
mathematically treatable, as particularly highlighted in Chen, Feldman \cite{ChenFeldmanpaper, ChenFeldman}, Elling, Liu \cite{EllingLiu, EllingLiu1}, 
see also the references listed in Dafermos \cite{Dafermos}, Ch. 18. 
In fact, this possibility remains limited, mostly, to flows that are irrotational, 
isentropic, for the polytropic gas law, 
and described in self-similar variables. 

\medskip
In particular, for irrotational flows, discontinuities of the vortex sheet 
type are not admissible. Since such discontinuities are the primary 
source of instabilities in compressible flows, assuming the flow is 
potential removes the possibility of accurately capturing the 
asymptotic behavior of supersonic into supersonic or harder supersonic into subsonic regimes.
Therefore, the study of 
systems of conservation laws, even in two spatial variables, is still in its infancy, and we bring in this paper 
some new directions in this context, for fluids which are non-potential (also called rotational). 

\medskip
Our development in Section \ref{SelfSimSysRF} relies on the Hodge-Helmholtz decomposition for vector fields in
connected domains, and typically we assume sufficient regularity of the functions involved. 
These assumptions are justified since our focus is on solutions of the Euler system in the 
subsonic regimes. More precisely, when the compressible 2-D Euler system is written in 
self-similar variables, particularly with the use of pseudo-velocities, 
the system is of mixed type. This means that, the system is strictly hyperbolic if, and only if,  
in pseudo-supersonic regimes. 
For real applications, such as the Prandtl–Meyer problem, the pseudo-supersonic regimes are governed 
by the Rankine–Hugoniot conditions. In contrast, in the pseudo-subsonic regimes one must handle a degenerate 
elliptic partial differential equation arising in the case of potential flow. For generalized polytropic gases, 
see \eqref{pressurelaw} below, 
such an equation does not satisfy a standard maximum principle (for $\gamma< 1$). 
The situation becomes even more intricate for rotational flows, where the 
governing equations form a coupled system. 

\medskip
In order to develop our study of the behavior of rotational flows, we define the well-known physical quantity 
called enthalpy, which is a strictly increasing function of the density. In this way, the density can be explicitly 
expressed in terms of the pseudo-velocity. To achieve this, we establish a Bernoulli type equation 
that allows the sound speed to be written in terms of the pseudo-velocity vector field.
Until now we have not fixed any
specific pressure law, which is necessary to write explicitly the sound speed in term of the pseudo-velocity. 
Since we are assuming that the process is isentropic, the pressure
depends only on the density, that is to say, a barotropic gas law. 
In this paper, we consider the pressure law $p=p(\rho)$, namely generalized polytropic gas, as follows
\begin{equation}
\label{pressurelaw}
p(\rho)=\dfrac{a^2}{\gamma}(\rho^{\gamma}- \underline{\rho}^{\gamma}),
\end{equation}
where $a>0$, $\gamma\in [-1, + \infty) \backslash  \{0\}$, $\rho > \underline{\rho} \ge 0$, $(\underline{\rho}> 0$ for $\gamma< 1)$. 
This pressure law generalizes the polytropic gas law, ($\gamma> 1$), and include the Chaplygin gas, ($\gamma= -1$), while 
ensuring that the pressure remains positive, even for negative values of gamma, which is a 
fluid dynamics condition,
we refer to Serre \cite{Serre2009}, p. 543. 
The equation above satisfies, for any $\gamma$, the important condition 
$p'(\rho)>0$. This is the most important condition from the pressure law in order that the 
results established here be true. 
Within this a slight modification of the pressure law considered, while ensuring that  $p'(\rho)> 0$,
the pressure is negative which allows for the inclusion of dark energy effects, that is, for $\gamma \in [-1, 0)$
\begin{equation}
\label{pressurelawmod}
    p(\rho)= - a^2 (\rho^{\gamma}- \underline{\rho}^{\gamma}) \equiv C_0 - C_1 \rho^{\gamma}. 
\end{equation}
Thus, by substituting the linear momentum equation into the continuity equation 
written in self-similar variables, we obtain a scalar equation that depends solely 
on the pseudo-velocity. 
To close the system, we observe that an additional transport equation for the vorticity 
has been derived, this constitutes the main ingredient used to derive system 
\eqref{newsystemU}, as stated in Theorem \ref{MainTheorem}.

\medskip
Section \ref{PotentialFlows} developed the theory for potential flows 
in the pseudo-subsonic regime, considering the generalized polytropic gas. 
In the case of potential flows, system \eqref{newsystemU} reduces to a scalar 
partial differential equation, which is elliptic precisely within the pseudo-subsonic regime.
First, we establish an Ellipticity Principle given by Theorem \ref{Ellipticity_principle_theorem}. 
One remarks that, this result was first established by Elling, Liu \cite{EllingLiu}
for $\gamma > -1$, and later extended by Chen, Feldman \cite{ChenFeldman},
who provided a more precise estimate for the case $\gamma \geq 1$. 
We were able to prove Theorem \ref{Ellipticity_principle_theorem}
for $\gamma \ge -1$ by closely following the approach of Chen-Feldman. 
Then we consider the Dirichlet problem \eqref{dirichlet}, which is a challenging 
free-boundary problem, since the domain $\Omega$, see \eqref{Omega}, depends on 
the solution itself. Moreover, the elliptic operator in \eqref{dirichlet} degenerates on 
the sonic line, that is, on $\partial \Omega$. 
Hence, we regularize it by considering the Dirichlet problem \eqref{qregularized}, 
whose solutions, for $\gamma \geq -1$, are established in Theorem~\ref{theorem4.1}. 
In particular, we remark that for $\gamma < 1$ we do not have a maximum principle,
therefore, refined techniques, such as the Moser iteration method combined 
with suitable estimates, are employed.
At the end of this section, we prove Theorem \ref{theorem4.3}, 
which establishes the existence of solutions to \eqref{dirichlet} 
for $\gamma > -1$. The proof relies fundamentally on the barrier 
provided by the Ellipticity Principle.

\medskip
Finally, we stress Section \ref{SecPM}, 
which is related to the analysis 
of the Euler equations for quasi-potential flows, a concept introduced in this section 
and also an important part of this paper. Indeed, a very special case occurs in situations 
where the amount of vorticity is small, more specifically, when the non-potential part of 
the pseudo-velocity vector field is significantly smaller than the potential part. 
In this context, we demonstrate in Theorem \ref{QPTHM} the existence of solutions to 
quasi-potential flows in the pseudo-subsonic regime. To achieve this, we employ a 
clever application of the Implicit Function Theorem on Banach spaces, where the 
Ellipticity Principle Theorem is once again fundamental in establishing a uniform 
ellipticity condition.
In particular, this section is interesting  
for numerical applications, where specific computationally efficient techniques 
already designed to deal with the potential flows of the Euler equations, 
can be employed to treat the case of quasi-potential flows.  

%%%%%%%%%%%%%%%%%%%%%%%
\subsection{2-D Isentropic Euler Equations}
\label{IEQ}
%%%%%%%%%%%%%%%%%%%%%%%

In this section, we recall some basic concepts of the theory 
of systems of conservation laws in two space dimensions. 
Moreover, we present the 2-D isentropic Euler equations.
In particular, we describe the case of steady flows 
at the end of this section. 
Most of the results are well-known, and we omit the proofs 
unless clarification is needed.

\medskip
To begin, a system of conservation laws in two space dimensions is a 
system of partial differential equations in the following form 
\begin{equation}
\label{2system}
U_t + \partial_x (F(U)) + \partial_y (G(U))=0,
\end{equation}
where $F,G:\R^m \to \R^m$ are given and $U: [0, \infty) \times \R^2 \to \R^m$ is the unknown, $U=U(t, x,y)$.
Writing equation (\ref{2system}) in the non-divergence form, we have
\begin{equation}
\label{2sysnondiv}
U_t+ DF(U) U_x + DG(U) U_y=0, 
\end{equation}
where we have assumed that, the functions involved are sufficiently smooth. 
Then, we may consider the following 
\begin{definition}
\label{defhyp2d}
The system \eqref{2system} is said to be hyperbolic if for every $U \in \R^m$ and each vector $\mathbf{\alpha}=(\alpha_1,\alpha_2) \in S^1$, the matrix
$$J(U , \mathbf{\alpha}) \equiv \alpha_1 DF(U) + \alpha_2 DG(U)$$
has $m$ real and distinct eigenvalues, where $DF, DG$ are the Jacobian matrices of the fluxes $F$ and $G$, respectively. 
Moreover, in this case, we write
\begin{equation}
\lambda_1(z)<\lambda_2(z) < \ldots < \lambda_m(z) \quad (z \in \R^m).
\end{equation}
\end{definition}

\medskip
Next we state a property that will be useful later to introduce the concept of self-similar variables in the context of the Euler equations.
\begin{proposition}\label{propinvariant}
The system \eqref{2system} is invariant under the self-similar scaling
$$
    (t,x) \mapsto (\mu t, \mu x), \quad \mbox{for any $ \mu \neq 0$}.
$$
\end{proposition}

\begin{comment}
%\begin{proof}
Consider $U(x,y,t)$ a solution of (\ref{2system}) and define
\begin{equation}
V(x,y,t):= U(\alpha x, \alpha y, \alpha t).
\end{equation}
Then,
\begin{eqnarray*}
\dfrac{\partial V}{\partial t}(x,y,t)&=& \dfrac{\partial U}{\partial t}(\alpha x, \alpha y, \alpha t) \alpha = \alpha U_t( \alpha x,\alpha y, \alpha t),\\
\dfrac{\partial V}{\partial x}(x,y,t)&=& \dfrac{\partial U}{\partial x}(\alpha x, \alpha y, \alpha t) \alpha = \alpha U_x( \alpha x,\alpha y, \alpha t),\\
\dfrac{\partial V}{\partial y}(x,y,t)&=& \dfrac{\partial U}{\partial y}(\alpha x, \alpha y, \alpha t) \alpha = \alpha U_y( \alpha x,\alpha y, \alpha t).
\end{eqnarray*}
So, we may write
\begin{equation*}
\begin{split}
& V_t + \partial_x (F(V)) + \partial y (G(V))\\[10pt]
&= V_t + DF(V) V_x + DG(V) V_y\\[10pt]
&= \alpha \left( U_t(\alpha x, \alpha y, \alpha t) + DF(U(\alpha x, \alpha y, \alpha t)) U_x(\alpha x, \alpha y, \alpha t)\right.\\[5pt]
& \quad + \left. DG(U(\alpha x, \alpha y, \alpha t))  U_y(\alpha x, \alpha y, \alpha t)\right)\\
&=0.
\end{split}
\end{equation*}
%\end{proof}
\end{comment} 

\medskip
Now, let us consider the isentropic Euler equations 
for compressible gas flow in a two-dimensional space, that is to say
\begin{equation}
\label{isentropic1}
\left\{
\begin{split}
& \partial_t \rho + \nabla_x \cdot (\rho \mathbf{v})= 0, \hspace{73pt} \mbox{(continuity equation)},
\\[5pt]
& \partial_t (\rho \mathbf{v}) + \nabla_x \cdot (\rho \mathbf{v} \otimes \mathbf{v}) + \nabla_x p= 0, \quad \mbox{(linear momentum equation)}.
\end{split}
\right.
\end{equation}
Here $t> 0$ is the time, $x=(x_1,x_2) \in \R^2$ 
is the space coordinates, $\rho=\rho(t,x)>0$ denotes the density, 
$\mathbf{v}=\mathbf{v}(t,x)=(v_1,v_2)$ is the velocity vector field, 
and $p=p(\rho) \in \R$ is the pressure law. 
Also, we recall that $\mathbf{a} \otimes \mathbf{b}$ denotes the tensor 
product of vectors $\mathbf{a}$ and $\mathbf{b}$. 
%that is, the matrix
%$$\mathbf{a} \otimes \mathbf{b} = [a_i b_j]_{ij} \quad (i,j=1,2).$$

\medskip
We will suppose $p'(\rho)>0$ 
and define $c:= \sqrt{p'(\rho)}$, which is called the sound speed. 
The flow is called subsonic, sonic or supersonic at $(t,x)$, respectively, if 
$$|\mathbf{v}(t,x)| < c(\rho(t,x)), \quad |\mathbf{v}(t,x)| = c(\rho(t,x))  \quad \mbox{or} \quad |\mathbf{v}(t,x)| > c(\rho(t,x)).$$
Moreover, 
we define the Mach number as
\begin{equation}
M:= \dfrac{|\mathbf{v}|}{c}.
\end{equation}
Then, the flow is subsonic, sonic or supersonic, respectively, if
$$
  M < 1, \quad M=1 \quad \mbox{or} \quad M>1.
$$

\medskip
In what follows, we examine the hyperbolicity of \eqref{isentropic1}.
Indeed, due to a straightforward calculation, the system \eqref{isentropic1} reads
%\begin{equation}
%\label{2deuler}
$$
\left\{
\begin{split}
& \rho_t + \partial_x(\rho u) + \partial_y (\rho v)= 0,
\\[5pt]
& (\rho u)_t + \partial_x(\rho u^2 + p(\rho))+ \partial_y (\rho u v)= 0,
\\[5pt]
& (\rho v)_t + \partial_x(\rho u v) + \partial_y(\rho v^2 + p(\rho))= 0,
\end{split}
\right.
%\end{equation}
$$
which can be rewritten in the form of equation \eqref{2system}. To this end, we denote
$$
\begin{aligned}
(U^1,U^2,U^3)&=(\rho, \rho u, \rho v), 
\\[5pt]
F(U^1,U^2,U^3)
%&=\left(U^2, \dfrac{(U^2)^2}{U^1} + p(U^1), \dfrac{U^2 \cdot U^3}{U^1}\right)\\
&=(\rho u, \rho u^2 + p(\rho), \rho u v)
\\[5pt]
G(U^1,U^2,U^3)
%&=&\left(U^3,\dfrac{U^2 \cdot U^3}{U^1}, \dfrac{(U^3)^2}{U^1} + p(U^1)\right)\\
&=(\rho v, \rho u v, \rho v^2+ p(\rho)).
\end{aligned}
$$
Taking $\alpha= (\alpha_1, \alpha_2) \in S^1$, we have
$J(U , \mathbf{\alpha})= \alpha_1 DF(U) + \alpha_2 DG(U)$, and 
calculating the eigenvalues, we obtain 
$$
\lambda_{1,3}=
u\alpha_1+v \alpha_2 \pm c, \quad 
\lambda_2= u\alpha_1 + v \alpha_2. 
$$
Since $c> 0$, it follows that 
$$\lambda_1 < \lambda_2 < \lambda_3,$$
from which we conclude that the system \eqref{isentropic1} is strictly hyperbolic.

\begin{remark}
Although we have verified that the system \eqref{isentropic1} is strictly hyperbolic, 
there are two types of physically relevant flows for which this property fails: 
the steady flow and the self-similar flow.

\smallskip
We will end this section by 
making a brief introduction to the former case, and the second one
will be detailed in the next section.
\end{remark}

Considering the Isentropic Euler equations for steady flows, it follows that, 
the density and the velocity vector field do not depend, explicitly, on time.
Thus, we have the system
\begin{equation}
\label{isentropicsteady}
\left\{
\begin{split}
&  \nabla_x \cdot (\rho \mathbf{v})= 0, 
\\[5pt]
& \nabla_x \cdot (\rho \mathbf{v} \otimes \mathbf{v}) + \nabla_x p= 0.
\end{split}
\right.
\end{equation}
Rewriting in the non-divergence form, we have
\begin{eqnarray}
\label{steadynondiv}
\left[ \begin{array}{ccc}
0&1&0\\
-u^2+ p'(\rho)&2u&0\\
-uv&v&u
\end{array}
\right]
\left[
\begin{array}{c}
\rho_x\\
(\rho u)_x\\
(\rho v)_x
\end{array}
\right]& & \nonumber\\
+\left[ \begin{array}{ccc}
0&0&1\\
-uv&v&u\\
-v^2+p'(\rho)&0&2v
\end{array}
\right]
\left[
\begin{array}{c}
\rho_y\\
(\rho u)_y\\
(\rho v)_y
\end{array}
\right]
&=&\left[\begin{array}{c}
0\\0\\0
\end{array}\right], 
\end{eqnarray}
that is to say, 
$A(U)U_x+B(U)U_y=0$,
and assuming that $A(U)$ is invertible, 
we have
$$
U_x+[A(U)]^{-1}B(U)U_y=0.
$$
One remarks that, this system can be interpreted as a one dimensional system of conservation laws,
and also it will be strictly hyperbolic if the eigenvalues of the matrix
$[A(U)]^{-1}B(U)$
are real and distinct.
Moreover, we note that 
\begin{eqnarray*}
\det (A^{-1}B-\lambda I)=0 &\Leftrightarrow & \frac{\det A \det (A^{-1}B- \lambda I)}{\det A}=0\\
& \Leftrightarrow & \det(B- \lambda A) =0.
\end{eqnarray*}
Thus, the eigenvalues are given by
$$
\begin{array}{rcl}
\lambda_2= \dfrac{v}{u} \quad \mbox{or} \quad  \lambda_{1,3} &=& \dfrac{-2uv \pm \sqrt{4u^2v^2 - 4 (c^2-u^2)(c^2-v^2)}}{2(c^2-u^2)}\\[10pt]
 %&=& \dfrac{-2uv\pm \sqrt{4u^2v^2-4c^4+4v^2c^2+4u^2c^2-4u^2v^2}}{2(c^2-u^2)}\\[10pt]
 &=& \dfrac{uv \pm c \sqrt{u^2+v^2-c^2}}{c^2-u^2}.
\end{array}
$$

\begin{remark}
\label{MixedEqs}
One remarks that, the system \eqref{isentropicsteady} is strictly hyperbolic provided $u^2+v^2 >c^2$, that is, in the supersonic regime.
Otherwise, the eigenvalues $\lambda_{1,3}$ are complex $u^2+v^2< c^2$, 
or equal $u^2+v^2= c^2$, respectively subsonic and sonic regimes. Consequently, the system \eqref{isentropicsteady} is of mixed type. 
\end{remark}

%%%%%%%%%%%%%%%%%%%%%%%%%%%%
\section{Self-similar System for Rotational Flows}
\label{SelfSimSysRF}
%%%%%%%%%%%%%%%%%%%%%%%%%%%%

The main issue of this section is to derive a
self-similar system for rotational flows, see \eqref{newsystemU}. 
To this end, we rewrite the Euler equations using self-similar variables, 
and apply the Hodge-Helmholtz decomposition 
for vector fields in connected domains. Indeed,
first, we show the existence of a 
Bernoulli type equation for flows which 
are not necessarily potential. Then, we consider 
a generalized polytropic gas law given by \eqref{pressurelaw}, obtain 
the vorticity equation and state the main theorem of this section.  

\medskip
Let us recall the isentropic Euler equations in two space dimensions given by system (\ref{isentropic1}), that is,
\begin{equation*}
\left\{
\begin{split}
& \partial_t \rho + \nabla_x \cdot (\rho \mathbf{v})= 0,
\\[5pt]
& \partial_t (\rho \mathbf{v}) + \nabla_x \cdot (\rho \mathbf{v} \otimes \mathbf{v}) + \nabla_x p= 0. 
\end{split}
\right.
\end{equation*}
Since first-order systems of conservations laws are invariant under the self-similar scaling (Proposition \ref{propinvariant}), 
we can seek for self-similar solutions, which depend only on the similarity coordinates. Conveniently, we take 
$$\mu= \frac{1}{t} \quad (t>0),$$
and denote
$$(\mathbf{v},p,\rho)(t,x_1,x_2)=(\mathbf{v},p,\rho)(\xi_1,\xi_2), \qquad \xi=(\xi_1,\xi_2)=\dfrac{x}{t}=\left(\dfrac{x_1}{t}, \dfrac{x_2}{t}\right).$$
This procedure is interesting in itself, as it reduces the number of variables.
Also, we define
\begin{equation}\label{defU}
U (\xi_1, \xi_2) \equiv (U^1,U^2)(\xi_1, \xi_2) := \mathbf{v}(\xi_1,\xi_2) - (\xi_1, \xi_2),
\end{equation}
called the pseudo-velocity. 

\medskip
Our goal now is to express the isentropic Euler equations in terms of the pseudo-velocity, 
which is typically done under the assumption of potential flow. 
Here, we do not adopt this approach and instead make the following important observation.

\begin{remark}[Hodge-Helmholtz decomposition] 
In most of the literature, the hypothesis that the pseudo-velocity $U$ is potential is used. 
That is, there exists $\psi=\psi(\xi)$, called potential, such that, $U=\nabla \psi$, where the gradient is taken 
with respect to the self-similar variables. 
In this case, we can also say that the flow is potential or the flow is irrotational. 

\smallskip
In this work, we employ the Hodge–Helmholtz decomposition to investigate the 
general case in which the pseudo-velocity $U$ is not necessarily potential. 
We emphasize that all results in this section remain valid if the flow is irrotational.
\end{remark}

Applying the Hodge-Helmholtz decomposition, there exist scalar and vector value functions, respectively,
$\psi(\xi) \in \R$, $W(\xi) \in \R^2$, $\xi= (\xi_1, \xi_2) \in \R^2$, such that, 
$$
   U= (U^1,U^2)= \nabla \psi + W, \quad \text{with  $\,\, \dive W= 0.$}
$$
Then, the flow will be rotational or irrotational (potential) if $W\neq 0$ or $W=0$, respectively.

\begin{remark}
When working with the pseudo-velocity, we will use 
superscript indices to indicate vector components, 
and subscript indices to indicate partial derivatives, that is,
$$\psi_i= \dfrac{\p \psi}{\p \xi_i}, \quad (i=1,2).$$
Therefore, we have 
$$U=(U^1,U^2)= \nabla \psi + W = (\psi_1, \psi_2) + (W^1,W^2), \qquad \dive W=0.$$
We also denote
\begin{equation}
U^{\perp}=(-U^2,U^1)
\end{equation}
and define 
\begin{equation}
\omega := \mathrm{rot}\, U = \dfrac{\p U^2}{\p \xi_1}- \dfrac{\p U^1}{\p \xi_2},
\end{equation}
where $\omega$ is called the vorticity of the fluid.
We notice that, under the Hodge-Helmholtz decomposition, 
$$\omega = \mathrm{rot}\, U = \mathrm{rot}\, W.$$
\end{remark}

\medskip
We say that the flow is pseudo-subsonic, pseudo-sonic or pseudo-supersonic at $(\xi_1, \xi_2)$ if
$$|U| < c(\rho(\xi_1, \xi_2)), \quad |U| = c(\rho(\xi_1, \xi_2))  \quad \mbox{or} \quad |U| > c(\rho(\xi_1, \xi_2)),$$
respectively.
Also, we define the pseudo-Mach-number as
\begin{equation}
L:= \dfrac{|U|}{c}.
\end{equation}
Then, the flow is pseudo-subsonic, pseudo-sonic or pseudo-supersonic at $(\xi_1, \xi_2)$ if
$$L < 1, \quad L=1 \quad \mbox{or} \quad L>1,$$
respectively.

\medskip
Now, we rewrite the isentropic Euler equations in terms of the pseudo-velocity, 
that is, we have the following 
\begin{proposition}
If $(\mathbf{v}, p, \rho)$ satisfy the system \eqref{isentropic1} 
and $U$ is defined by \eqref{defU} then, assuming enough regularity, it follows that
\begin{equation}
\label{selfisentropic}
\left\{
\begin{split}
&\dive (\rho \, U)+ 2 \, \rho=0,
\\[5pt] 
&\dive (\rho \, U \otimes U)+ Dp + 3 \, \rho \, U= 0.
\end{split}
\right.
\end{equation}
\end{proposition}

\begin{proof}
1. First, we note that

$$\dfrac{\p \xi_i}{\partial t}= - \xi_i \dfrac{1}{t}, \, \quad \dfrac{\p \xi_i}{\p x_j}=\delta_{ij}\dfrac{1}{t} \qquad(i,j=1,2),$$
where $\delta_{ij}$ is the Kronecker delta.
Then, it follows from the continuity equation $(\ref{isentropic1})_1$ that, 
%$$
 %  \partial_t \rho + \nabla_x \cdot (\rho \mathbf{v})= 0,
%$$
\begin{eqnarray*}
0&=& \p_t (\rho(\xi)) + \nabla_x \cdot (\rho(\xi) \mathbf{v}(\xi))\\
%&=& \dfrac{\p \rho}{\p \xi_1} \dfrac{\p \xi_1}{\p t}+ \dfrac{\p \rho}{\p \xi_2} \dfrac{\p \xi_2}{\p t}+ \dfrac{\p}{\p \xi_1}(\rho v_1) \dfrac{\p \xi_1}{\p x_1}+ \dfrac{\p}{\p \xi_2}(\rho v_2) \dfrac{\p \xi_2}{\p x_2}\\
%&=& \dfrac{\p \rho}{\p \xi_1}(-\xi_1)\dfrac{1}{t}+ \dfrac{\p \rho}{\p \xi_2}(- \xi_2)\dfrac{1}{t}+ \rho \dfrac{\p v_1}{\p \xi_1} \dfrac{1}{t}+ v_1 \dfrac{\p \rho}{\p \xi_1} \dfrac{1}{t}+ \rho \dfrac{\p v_2}{\p \xi_2} \dfrac{1}{t}+ v_2 \dfrac{\p \rho}{\p \xi_2} \dfrac{1}{t}\\
&=& - \nabla \rho \cdot \xi + \rho \dive \mathbf{v} + \mathbf{v} \cdot \nabla \rho \\
&=& - \xi \cdot \nabla \rho + \nabla \cdot (\rho \mathbf{v}).  
\end{eqnarray*}
Using $\mathbf{v}=U+\xi$, we have
$$
\begin{aligned}
0&= - \xi \cdot \nabla \rho + \dive (\rho U) + \dive (\rho \xi)\\
%&= - \xi \cdot \nabla \rho + \dive (\rho U) + \rho \dive \xi + \nabla \rho \cdot \xi\\
&= \dive (\rho U)+ 2 \rho.
\end{aligned}
$$

2. Now, fixing $i$ in the linear momentum equation, we have
$$
\begin{aligned}
0&= \p_t (\rho v_i)+ \nabla_x \cdot (\rho v_i \mathbf{v})+ \dfrac{\p p}{\p x_i}
\\[5pt]
&= \dfrac{\p}{\p \xi_1} ( \rho v_i) \dfrac{\p \xi_1}{\p t} + \dfrac{\p}{\p \xi_2}(\rho v_i) \dfrac{\p \xi_2}{\p t} + \dfrac{\p}{\p \xi_1}(\rho v_i v_1) \dfrac{\p \xi_1}{\p x}+ \dfrac{\p}{\p \xi_2}(\rho v_i v_2) \dfrac{\p \xi_2}{\p x_2}+ \dfrac{\p p}{\p \xi_i} \dfrac{\p \xi_i}{\p x_j}
%\\[5pt]
%&= \Big( \rho \dfrac{\p v_i}{\p \xi_1}+ v_i \dfrac{\p \rho}{\p \xi_1} \Big) (-\xi_1)\dfrac{1}{t}+ \Big( \rho \dfrac{\p v_i}{\p \xi_2}+ v_i \dfrac{\p \rho}{\p \xi_2} \Big) (-\xi_2) \dfrac{1}{t}
%\\[5pt]
%&+ \! \Big( \rho v_i \dfrac{\p v_1}{\p \xi_1}+ v_1 \rho \dfrac{\p v_i}{\p \xi_1}+ v_1 v_i \dfrac{\p \rho}{\p \xi_1} \Big) \dfrac{1}{t} \!
%+\!  \Big( \rho v_i \dfrac{\p v_2}{\p \xi_2}+ v_2 \rho \dfrac{\p v_i}{\p \xi_2}+ v_2 v_i \dfrac{\p \rho}{\p \xi_2} \Big) \dfrac{1}{t} \! + \dfrac{\p p}{\p \xi_i}\dfrac{1}{t}
\\[5pt]
&= \rho (-\xi) \cdot \nabla v_i + v_i (- \xi) \cdot \nabla \rho + \rho v_i \dive \mathbf{v} + \rho (\mathbf{v} \cdot \nabla v_i) + v_i(\mathbf{v} \cdot \nabla \rho) + \dfrac{\p p}{\p \xi_i}.
\end{aligned}
$$
Therefore, we obain 
\begin{equation*}
\begin{split}
&\overbrace{\rho U^i \dive U + \rho (U \cdot \nabla U^i)+U^i (U \cdot \nabla \rho)}^{= \dive (\rho U^i U)}
\\
& + \dfrac{\p p}{\p \xi_i}+ 3 \rho U^i + \xi_i \underbrace{\left(\rho \dive U +U \cdot \nabla \rho + 2 \rho \right)}_{=0}=0
%\\
%& \dive (\rho U^i U) + \dfrac{\p p}{\p \xi_i}+ 3 \rho U^i=0, \qquad (i=1,2), 
\end{split}
\end{equation*}
that is to say, $\dive (\rho U \otimes U) +Dp + 3 \rho U= 0$. 
\end{proof}

We notice that, rewriting the system \eqref{selfisentropic} as 
$$
\left\{
\begin{split}
&\dive (\rho \, U)= -2 \rho,
\\[5pt]
&\dive (\rho \, U \otimes U)+ \nabla p= - 3 \rho \, U, 
\end{split}
\right.
$$
the terms of higher order are the same of the system \eqref{isentropicsteady}. 
Therefore, replacing $|U|^2$ by $L^2c^2$ the eigenvalues are given by
$$
\lambda_{1,3} = \dfrac{U^1U^2 \pm c^2 \sqrt{L^2-1}}{c^2-(U^1)^2}, \quad \lambda_2= \dfrac{U^2}{U^1}. 
$$
Consequently, the system \eqref{selfisentropic}, (similar to system \eqref{isentropicsteady}), is of mixed type. 
More precisely, this system is strictly hyperbolic, 
if and only if $L^2 > 1$, that is, in the pseudo-supersonic regime.
Otherwise, the eigenvalues $\lambda_{1,3}$ are complex $L^2< 1$, 
or equal $L^2= 1$, respectively pseudo-subsonic and pseudo-sonic regimes. 

\medskip
Now we claim that, if $U$ satisfies the above system, 
then substituting the linear momentum equation into the continuity equation, 
we obtain the following 
\begin{proposition}
If $U$ satisfies the system \eqref{selfisentropic}, then $U$ satisfies the equation
\begin{equation}
\label{eqgeralU}
c^2 \, \dive U - (DU)\, U \cdot U = |U|^2- 2 \, c^2.
\end{equation}
\end{proposition}

\begin{proof}
1. First, fixing $i$ in $(\ref{selfisentropic})_2$, we have
\begin{eqnarray*}
0&=& \dive (\rho U^i U)+ \dfrac{\p p}{\p \xi_i}+ 3 \rho U^i\\
&=&\rho U^i \dive U + \rho (U \cdot \nabla U^i)+U^i(U \cdot \nabla \rho) + \dfrac{\p p}{\p \xi_i} + 3 \rho U^i\\
%&=& U^i(\rho \dive U + U \cdot \nabla \rho + 2 \rho) + \rho (U \cdot \nabla U_i)+ \dfrac{\p p}{\p \xi_i}+ \rho U^i\\
&=& U^i\underbrace{(\dive (\rho U) + 2 \rho)}_{=0} + \rho (U \cdot \nabla U^i)+ \dfrac{\p p}{\p \xi_i}+ \rho U^i\\
%&=& \rho (U \cdot \nabla U^i)+ \dfrac{\p p}{\p \xi_i}+ \rho U^i\\
&=& \rho (U \cdot \nabla U^i)+ p'(\rho) \dfrac{\p \rho}{\p \xi_i} + \rho U^i.
\end{eqnarray*}
Then, we obtain, 
$\rho \, (DU \, U) + p'(\rho) \nabla \rho + \rho \,U= 0$, that is,
\begin{equation}
\label{eq1}
p'(\rho) \nabla \rho = - \rho \, U - \rho \, (DU\, U).
\end{equation}

\medskip
2. Now multiplying, $\dive (\rho \, U) + 2 \, \rho= 0$
by $p'(\rho)=c^2$, and using $$\dive (\rho U)= \rho \dive U + U \cdot \nabla \rho,$$ it follows that 
\begin{equation}\label{eq2}
c^2 \rho \, \dive U + p'(\rho) \nabla \rho \cdot U + 2 \, \rho \, c^2= 0.
\end{equation}
Substituting \eqref{eq1} into \eqref{eq2}, we have 
$$ 
   c^2 \rho \dive U + (- \rho U - \rho (DU) U) \cdot U + 2 \rho c^2= 0,
$$
that is,    
$$ 
   c^2 \, \rho \, \dive U - \rho \,|U|^2 - \rho \, (DU) \, U \cdot U + 2 \, \rho \, c^2=0.
$$
Finally, dividing by $\rho> 0$ we obtain \eqref{eqgeralU}.
\end{proof}

\medskip
We remark that, the above result, equation  \eqref{eqgeralU},  
is valid in the general case where $U= \nabla \psi + W$, with $\dive \, W= 0$.
In particular, when $W \equiv 0$ (irrotational flow), 
we recover the case where the pseudo-velocity is potential, 
and the equation \eqref{eqgeralU} takes the following well-known form,
\begin{equation}
\label{edpp}
c^2 \Delta \psi - (D^2 \psi) \nabla \psi \cdot \nabla \psi = |\nabla \psi|^2 - 2c^2,
\end{equation}
see for instance Dafermos \cite{Dafermos}, Chapter 18.

\medskip
Our goal now is to show that $c^2=c^2(\rho)$ can be written as $c^2=\mathcal{F}(U)$. 
Hence the equation \eqref{eqgeralU} is an equation depending only on $U$.
In order to show that, first we derive a Bernoulli type Equation, and then 
we consider a specific pressure law.   

%%%%%%%%%%%%%%%%%%
\subsection{Bernoulli type Equation}
\label{BernTypeEq}
%%%%%%%%%%%%%%%%%%
 
To begin, let us recall the definition of the important 
themodynamic quantity called enthalpy, $h= h(\rho)$, (isentropic process), given by
\begin{equation}
h'(\rho)=\dfrac{p'(\rho)}{\rho}>0.
\end{equation}
Note that, by definition, $h(\cdot)$ is strictly 
increasing and, therefore, invertible.
Then, we use the enthalpy to obtain a very important equation, 
called here a Bernoulli type equation.
To achieve this, we will employ the following three lemmas.
\begin{lemma}
\label{lemmaid1}
Let $R,S:\R^2 \to \R^2$ be smooth functions, denoted by $R= (R^1,R^2)$, and $S= (S^1,S^2)$. Then, 
$$
   \nabla (R \cdot S)= (DR)^TS +(DS)^T R.
$$
\end{lemma}
\begin{proof}
%First, we define
%$$
%    R_i := \left( R^1_i, R^2_i \right) = \left( \dfrac{\partial R^1}{\partial \xi_i}, \dfrac{\partial R^2}{\partial \xi_i} \right), 
%    \quad S_i:= \left( S^1_i, S^2_i \right) = \left( \dfrac{\partial S^1}{\partial \xi_i}, \dfrac{\partial S^2}{\partial \xi_i}\right), \quad (i=1,2). 
%$$
A direct calculation yields, 
\begin{eqnarray*}
\nabla (R \cdot S) 
%&=& \nabla (R^1S^1+ R^2S^2)\\
&=& (R^1S^1_1 + S^1R^1_1 + R^2S^2_1 +S^2R^2_1, R^1S^1_2 +S^1R^1_2+R^2S^2_2+S^2R^2_2)\\
%&=&(R_1 \cdot S + S_1 \cdot R, R_2 \cdot S + S_2 \cdot R)\\
%&=& (R_1 \cdot S, R_2 \cdot S)+(S_1 \cdot R, S_2 \cdot R)\\
&=& (DR)^T S + (DS)^T R.
\end{eqnarray*}
\end{proof}
\begin{lemma}\label{lemmaid2}
If $U=\nabla \psi + W$ is given by the Hodge-Helmholtz decomposition, then,
\begin{eqnarray}
& \mathrm{(1)}& \nabla ( \nabla \psi \cdot W)=(D^2 \psi) W+ (DW)^T \nabla \psi,\\
& \mathrm{(2)}& \nabla \left( \frac{1}{2}|W|^2 \right) =(DW)^T W,\\
& \mathrm{(3)}& \nabla \left( \frac{1}{2}|\nabla \psi|^2 \right)=(D^2 \psi)\nabla \psi,\\
& \mathrm{(4)}& (DW)U- (DW)^T U= \omega U^{\perp}.
\end{eqnarray}
\end{lemma}
\begin{proof}
The first three identities are directly derived from the previous lemma,
so we prove the last one. 
%Here we recall that, for $i = 1, 2$, 
%$$
 %   \nabla W^i= \left( W^i_1, W^i_2 \right) = \left( \dfrac{\partial W^i}{\partial \xi_1}, \dfrac{\partial W^i}{\partial \xi_2} \right), 
  %  \quad W_i= \left( W^1_i, W^2_i \right) = \left( \dfrac{\partial W^1}{\partial \xi_i}, \dfrac{\partial W^2}{\partial \xi_i} \right).
%$$
%Then, we obtain by 
Due to a simple algebraic manipulation, we obtain 
$$
\begin{aligned}
& (DW)U-(DW)^TU\\
&=  ( \nabla W^1 \cdot U - W_1 \cdot U, \nabla W^2 \cdot U - W_2 \cdot U )
\\
&= ( W^1_1 U^1+W^1_2U^2-W^1_1U^1 - W^2_1U^2, W^2_1U^1+W^2_2U^2-W^1_2U^1-W^2_2U^2)
\\
&= (- \omega U^2, \omega U^1)= \omega U^{\perp}.
\end{aligned}
$$
\end{proof}
\begin{lemma}\label{lemmaF}
If $F: \R^2 \to \R$ is a $C^2$ function and $\nabla F$ is given by
$$
  \big(\nabla F\big)(\xi_1,\xi_2)= \big(G(\xi_1,\xi_2), H(\xi_1,\xi_2)\big),
$$
then $F$ can be written as a function of $G,H$ by the following formula
\begin{equation}\label{FdepU}
F(\xi_1,\xi_2)= \int_0^{\xi_1} G(s,0) ds + \int_0^{\xi_2} H(\xi_1,s)ds+ C,
\end{equation}
where $C \in \R$ is a constant. 
\end{lemma}
\begin{proof}
Fixed $(\xi_1,\xi_2)\in \R^2$ and considering the curves $C_1$, $C_2$ parametrized by
$$
 C_1: \, \sigma_1(s)=(s,0), \, 0 \le s \le \xi_1, \quad 
 C_2: \, \sigma_2(s)=(\xi_1,s), \, 0 \le s \le \xi_2,
$$
respectively, we have
$$
\begin{aligned}
\int_0^{\xi_1} G(s,0)ds &+ \int_0^{\xi_2} H(\xi_1,s)ds
\\
%&=&\int_0^{\xi_1} (G(s,0),H(s,0))\cdot (1,0) ds + \int_0^{\xi_2} (G(\xi_1,s),H(\xi_1,s))\cdot (0,1) ds\\[5pt]
&= \int_0^{\xi_1} \nabla F (\sigma_1(s))\cdot \sigma'_1(s)ds + \int_0^{\xi_2} \nabla F(\sigma_2(s))\cdot \sigma'_2(s)ds
\\
%&=&\int_{(0,0)}^{(\xi_1,0)} \nabla F \cdot dr + \int_{(\xi_1,0)}^{(\xi_1,\xi_2)} \nabla F \cdot dr\\[5pt]
&= \int_{(0,0)}^{(\xi_1,\xi_2)} \nabla F \cdot dr= F(\xi_1,\xi_2)+C, \quad (C \in \R).
\end{aligned}
$$
\end{proof}

Now we consider the main theorem of this section.
\begin{theorem}
(Bernoulli type equation) If $U=\nabla \psi + W$ is a $C^2$ function that satisfy 
the system \eqref{selfisentropic}, then the following equation is satisfied
\begin{equation}
\label{bernoullistypeequation}
h+\psi+\frac{1}{2}|U|^2= \int_0^{\xi_1} G(s,0) ds + \int_0^{\xi_2} H(\xi_1,s)ds + C, \quad C \in \R,
\end{equation}
where
$$
  G:= \omega U^2-W^1, \quad H:= \omega U^1-W^2.
$$
\end{theorem}

\begin{proof}
1. First, returning to the equation (\ref{eq1}) and using the enthalpy $h$, the Hodge-Helmholtz decomposition and Lemmas \ref{lemmaid1} and \ref{lemmaid2}, we have
\begin{eqnarray*}
0&=& p'(\rho) \nabla \rho + \rho U + \rho (DU)U\\[5pt]
%&=& h'(\rho) \nabla \rho + U + (DU)U\\[5pt]
%&=& \nabla h + U +(DU)U\\[5pt]
&=& \nabla h + \nabla \psi + W + (D^2 \psi) \nabla \psi + (D^2 \psi) W + (DW)U\\[5pt]
%&=& \nabla h + \nabla \psi + (D^2 \psi) \nabla \psi + (D^2\psi)W + (DW)^T \nabla \psi +
%(DW)^TW+W\\
&&+(DW)U - (DW)^T\nabla \psi - (DW)^T W\\[5pt]
&=& \nabla h + \nabla \psi + \nabla \left( \frac{1}{2}|\nabla \psi|^2 \right) + \nabla (\nabla \psi \cdot W)+ \nabla \left( \frac{1}{2} |W|^2 \right) +W\\
&&+(DW)U - (DW)^T U\\[5pt]
%&=& \nabla \left( h+ \psi + \frac{1}{2}|\nabla \psi|^2 + \nabla \psi \cdot W + \frac{1}{2} |W|^2 \right) + W + \omega U^{\perp}\\[5pt]
&=& \nabla \left( h + \psi + \dfrac{1}{2}|U|^2 \right) + W + \omega U^{\perp}.
\end{eqnarray*}
Conveniently, we define
\begin{equation}\label{eqdeff}
F:=h + \psi + \frac{1}{2}|U|^2, 
\end{equation}
hence we obtain 
\begin{equation}\label{eqgradf}
\nabla F + W + \omega U^{\perp}=0.
\end{equation}

2. Now, using the above equation we can write
$\nabla F=(G,H)$, 
where
$$
G:= \omega U^2-W^1, \quad 
H:= \omega U^1-W^2, 
$$
which are functions of $U$. 
Then, using Lemma \ref{lemmaF}, the function $F$ can be written as a function depending only on $U$ by the formula
\begin{equation}\label{FdepU}
F(\xi_1,\xi_2)= \int_0^{\xi_1} G(s,0) ds + \int_0^{\xi_2} H(\xi_1,s)ds+ C, \quad C \in \R.
\end{equation}
Therefore, it follows from equations \eqref{eqdeff} and \eqref{FdepU} that
$$
    h+\psi+\frac{1}{2}|U|^2= \int_0^{\xi_1} G(s,0) ds + \int_0^{\xi_2} H(\xi_1,s)ds + C, \quad C \in \R,
$$
which completes the proof.
\end{proof}

\begin{remark}
We notice that if $W=0$, then we have $U=\nabla \psi$, that is, the fluid is potential and equation \eqref{eqgradf} turns out
$\nabla F=0$. Therefore, the equation \eqref{eqdeff} assumes the form
\begin{equation}
\label{bernoulli}
h+\psi+\dfrac{1}{2}|\nabla \psi|^2=C, \quad C \in \R, 
\end{equation}
which is the well-known Bernoulli's equation, see Dafermos \cite{Dafermos}, Chapter 18.
Therefore, we call equation \eqref{bernoullistypeequation} Bernoulli-type equation.  
\end{remark}

We finish this section observing that,  
since $h$ is an invertible function of $\rho$, it follows from the Bernoulli type equation that
\begin{equation*}
h(\rho)= \int_0^{\xi_1} G(s,0) ds + \int_0^{\xi_2} H(\xi_1,s)ds-\psi-\frac{1}{2}|U|^2 + C, \quad C \in \R,
\end{equation*}
and thus
\begin{equation}\label{rhoh}
\rho = h^{-1} \Big( \int_0^{\xi_1} G(s,0)ds + \int_0^{\xi_2} H(\xi_1,s)ds - \psi - \frac{1}{2}|U|^2 +C \Big). 
\end{equation}

%%%%%%%%%%%%%%%%%%%%%%%%%
\subsection{Generalized polytropic gas law}
%%%%%%%%%%%%%%%%%%%%%%%%%

As previously mentioned, we will use the Bernoulli type equation in order 
to write $c^2=\mathcal{F}(U)$, so the equation (\ref{eqgeralU}) is an 
equation depending only on $U$ instead of depending on $(\rho,U)$. 
We recall that, until now, we have not fixed any specific pressure law. 
This is necessary now in order to find the function $\mathcal{F}$.

\medskip
From now on we will fix the pressure law $p= p(\rho)$ as given by \eqref{pressurelaw}, that is, 
$$
p(\rho)=\dfrac{a^2}{\gamma}(\rho^{\gamma}- \underline{\rho}^{\gamma}),
$$
where $a>0$, $\gamma\in [-1, + \infty) \backslash  \{0\}$, $\rho > \underline{\rho} \ge 0$, $(\underline{\rho}> 0$ for $\gamma< 1)$.
%This pressure law generalizes the polytropic gas law and include the Chaplygin gas, that is $\gamma= -1$, while 
%ensuring that the pressure remains positive, even for negative values of gamma, which is a 
%fluid dynamics condition.
%We refer to Serre \cite{Serre2009}, p. 543. 
%For any $\gamma$, the important condition 
%$p'(\rho)>0$ is satisfied. 

\begin{remark}
One observes that, the equation \eqref{pressurelaw} satisfies, for each $\gamma$, the condition 
$p'(\rho)>0$. 
%This is the most important condition from the pressure law in order that the results established here be true. 
In fact, 
with a slight modification of \eqref{pressurelaw}, see \eqref{pressurelawmod} in the introduction,  
it is possible for the pressure to be negative for $\gamma \in [-1, 0)$.
This modification of the constitutive equation enables the incorporation of dark matter.
Indeed, the parameter $\gamma$
allows flexibility in the model, controlling the transition between matter-like and dark energy-like behavior. Lower values of 
$\gamma$ delay the transition to dark energy-like behavior.
\end{remark}

\begin{proposition}
If $(\mathbf{v},p,\rho)$ satisfy the system \eqref{isentropic1}, where 
$\mathbf{v}$ is $C^2$, $p$ is given by \eqref{pressurelaw}, 
$c$ is the sound speed and $U$ is the pseudo-velocity, 
then there exists a $C^2$ function $\mathcal{F}$, such that
\begin{equation}
c^2= \mathcal{F}(U).
\end{equation}
\end{proposition}
\begin{proof}
If $\gamma \neq 1$, then we have 
\begin{equation}\label{c21}
c^2(\rho)=p'(\rho)=a^2 \rho^{\gamma-1}, 
\end{equation}
and thus 
\begin{equation}
\int_{\underline{\rho}}^{\rho} \dfrac{p'(z)}{z} dz
= \int_{\underline{\rho}}^{\rho} a^2 z^{\gamma-2} dz
= h\left(\rho\right) - h\left(\underline{\rho}\right),
\end{equation}
where for each $z\geq 0$,
$$
     h\left(z\right)= \dfrac{a^2 z^{\gamma-1}}{\gamma-1}\geq 0.
 $$
Substituting the above equation into \eqref{c21}, we have 
\begin{equation}
\label{cfunctionh}
    c^2(\rho)= (\gamma-1) \left( h(\rho) - h\left(\underline{\rho}\right)\right)> 0. 
\end{equation}
Therefore, it follows from equations (\ref{cfunctionh}) and (\ref{rhoh}) that, 
$c^2=\mathcal{F}(U)$. 

\medskip
Now, if $\gamma=1$, then
\begin{equation}
p(\gamma)=a^2(\rho - \underline{\rho})
\end{equation}
and we have 
\begin{equation}
c^2(\rho)=p'(\rho)=a^2.
\end{equation}
Moreover, we note that in this case, since $c^2$ is a constant function 
the Bernoulli type equation is not required to show that $c^2=\mathcal{F}(U)$.
\end{proof}

%%%%%%%%%%%%%%%%%%%%%%%%%%%%
\subsection{The vorticity equation and main theorem}
%%%%%%%%%%%%%%%%%%%%%%%%%%%%

We have already proved that,
the equation \eqref{eqgeralU} depends only on $U=(U^1, U^2)$. 
Now, we will see that $U$ must satisfy an important condition 
called the vorticity equation, which we state below.
\begin{proposition}
If $U=\nabla \psi + W$ is a $C^2$ function that satisfy system \eqref{selfisentropic}, 
then the function $U$ must satisfy the vorticity equation, which is given by
\begin{equation}
\label{conditionU}
\dive (\omega U) + \omega =0.
\end{equation}
\end{proposition}

\begin{proof}
We have seen in the proof of the Theorem \ref{bernoullistypeequation} that there exists a $C^2$ function $F$ such that
$$\nabla F=(G,H),$$
where $G=\omega U^2 - W^1$ and $H=\omega U^1 - W^2$.
In order for such a function to exist, it is necessary that
$$ \dfrac{\p H}{\p \xi_1} = \dfrac{\p G}{\p \xi_2}.$$
Therefore, $U$ must satisfy the following condition
$$
 - \omega (\psi_{11} + W^1_1) - (\psi_1+W^1) \omega_1 - W^2_1=
\omega (\psi_{22}+W^2_2)+ (\psi_2+W^2) \omega_2 - W^1_2, 
$$
that is to say, 
$$ 
 \omega ( \underbrace{\psi_{11}+\psi_{22}}_{= \Delta \psi} + \underbrace{W^1_1+W^2_2}_{= \dive W=0} ) +
\underbrace{(\psi_1+W^1)\omega_1 + (\psi_2 + W^2) \omega_2}_{= U \cdot \nabla \omega} 
+ \underbrace{W^2_1 - W^1_2}_{= \omega}= 0.
$$
Consequently, we obtain,
$\omega (\Delta \psi + 1) + U \cdot \nabla \omega= 0$, or 
$$
   \dive (\omega U) + \omega= 0,
$$
which finish the proof. 
\end{proof}

\begin{remark}
1. In the context of the full Euler equations, where $p= p(\rho, s)$ with $s$ the entropy, instead of equation \eqref{conditionU}, we have
$$  \omega (\Delta \psi + 1) + U \cdot \nabla \omega = - \curl \left( \dfrac{Dp}{\rho} \right), $$
see Chen, Feldman \cite{ChenFeldman}, Chapter 18.

\medskip
2. If the flow is potential, 
then equation \eqref{conditionU} is trivially satisfied.
\end{remark}

Therefore, we have proved in this section the following.
\begin{theorem}
\label{MainTheorem}
Suppose that $(\mathbf{v},p,\rho)$ satisfies the isentropic Euler equations \eqref{isentropic1} 
and that $U$, the pseudo-velocity, is a $C^2$ function. Suppose also that $c^2=p'(\rho)> 0$ 
is the sound velocity, $h=h(\rho)$ is the enthalpy and the pressure law $p=p(\rho)$ is given by
$$p(\rho)=\dfrac{a^2}{\gamma}(\rho^{\gamma} - \underline{\rho}^{\gamma}),$$
where $a>0$, $\gamma\in [-1, + \infty) \backslash  \{0\}$, $\rho > \underline{\rho} \ge 0$, $(\underline{\rho}> 0$ for $\gamma< 1)$.
\begin{enumerate}
\item Then, there exists $\mathcal{F}$ such that $c^2=\mathcal{F}(U)$ and the pseudo-velocity $U$ satisfies the system
\begin{equation}
\label{newsystemU}
\left\{
\begin{split}
& c^2 \dive U - (DU) U \cdot U = |U|^2 - 2c^2,
\\[5pt]
& \dive (\omega U) + \omega =0,
\end{split}
\right.
\end{equation}
where $\omega = \mathrm{rot}\, U$.

\item Applying the Hodge-Helmholtz decomposition, that is,  
$U= \nabla \psi + W$,
with $\dive W=0$, there exists $F$, such that,
$\nabla F = - \omega U^{\perp} - W$. 
Moreover, $c^2(U)$ is given by
\begin{equation}\label{c2U}
c^2(U)=(\gamma-1) \left( F - \psi - \dfrac{1}{2}|U|^2 \right), \quad \text{with $C= 0$},
\end{equation}

if $\gamma \neq 1$, and
$c^2=a^2,$ if $\gamma=1$.
Also, holds the Bernoulli type equation
\begin{equation}
h+ \psi + \dfrac{1}{2}|U|^2 = F, \quad \text{with $C= 0$}.
\end{equation}
\end{enumerate}
\end{theorem}

\begin{definition}
Under the conditions of Theorem \ref{MainTheorem}, system \eqref{newsystemU}
is called the self-similar system for rotational flows.  
\end{definition}

%%%%%%%%%%%%%%%%%%%%%%
\section{Potential Flows for a Generalized Polytropic Gas}
\label{PotentialFlows}
%%%%%%%%%%%%%%%%%%%%%%

In this section, we assume 
$W= 0$ in the Hodge–Helmholtz decomposition and consider the potential case 
$U= \nabla \varphi$, under the assumption that the pressure law is given by \eqref{pressurelaw}. 
We will focus on the pseudo-subsonic regimes, where the equation \eqref{edpp} is elliptic as described below. In this case
we remark that, for $-1 \leq \gamma < 1$, the associated second-order (elliptic) partial differential equation involves a zeroth-order term 
with positive sign, so that the maximum principle may fail. 
This feature renders the problem of establishing the existence 
of solutions more challenging than in the classical polytropic regime $\gamma \geq 1$.

\medskip 
In the potential setting, systems \eqref{selfisentropic} reduces to
\begin{equation}
\left\{
\begin{array}{rcl}
\dive (\rho \nabla \varphi) + 2 \rho &=& 0,
\\[5pt]
\dive (\rho \nabla \varphi \otimes \nabla \varphi ) + Dp + 3 \rho \nabla \varphi &=&0, 
\end{array}\right.
\end{equation}
and we recall the equation \eqref{edpp}, that is, 
$$
%\label{edpp}
c^2 \Delta \varphi - \sum_{i,j=1}^2 \varphi_i \varphi_j \varphi_{ij} = | \nabla \varphi|^2 - 2c^2,
$$
where
\begin{equation}
\label{eqc2}
\begin{aligned}
c^2(\varphi)&= - (\gamma-1)\left( \varphi + \frac{1}{2}|\nabla \varphi|^2 \right), \quad \text{if $\gamma \neq 1$}, 
\\[5pt] 
c^2(\varphi)&= a^2, \quad \text{if $\gamma=1$.}
\end{aligned}
\end{equation}
Moreover, 
we can rewrite equation \eqref{edpp} in the form
$$A \varphi_{11} + B \varphi_{12} + C \varphi_{22} + \; \mbox{lower order terms}= 0,$$
in order to classify it as an elliptic, parabolic or hyperbolic second order partial differential equation. Indeed, we have
$$ 
    \left( 1 - \dfrac{\varphi_1^2}{c^2} \right) \varphi_{11}  -2 \dfrac{\varphi_1 \varphi_2}{c^2} \varphi_{12} + \left(1 - \dfrac{\varphi_2^2}{c^2} \right) \varphi_{22} - \dfrac{|\nabla \varphi |^2}{c^2}+2= 0.
$$
Recall that,
$L^2 \; c^2 = |\nabla \varphi|^2$, 
and therefore we have
\begin{eqnarray*}
B^2-4AC &=& \left( -2 \dfrac{\varphi_1 \varphi_2 L^2}{|\nabla \varphi|^2} \right)^2 - 4 \left( 1 - \dfrac{L^2 \varphi_1^2}{|\nabla \varphi|^2} \right) \left( 1 - \dfrac{L^2 \varphi_2^2}{|\nabla \varphi|^2} \right)\\
&=& \dfrac{1}{|\nabla \varphi|^4} \left( 4 \varphi_1^2 \varphi_2^2 L^4 - 4 (|\nabla \varphi|^2 - L^2 \varphi_1^2)(|\nabla \varphi|^2 - L^2 \varphi_2^2) \right)
\\
&=& \dfrac{4}{|\nabla \varphi|^2} \left( L^2 \varphi_1^2 +L^2\varphi_2^2 - | \nabla \varphi|^2 \right)
%\\
= 4(L^2-1). 
\end{eqnarray*}
Then, equation \eqref{edpp} is elliptic, parabolic, or hyperbolic according as $L^2 < 1$, $L^2 = 1$, or $L^2 > 1$, respectively.
It follows that, in the case of potential flow, the classification of \eqref{edpp} 
aligns more precisely with that of system \eqref{selfisentropic}. 

\medskip
Together with equation \eqref{edpp} considered in the pseudo-subsonic regimes, we assume the Dirichlet boundary condition, namely,
$$
   \varphi = \varphi_b \quad \mbox{on} \quad \partial \Omega,
$$
where $\vp_b$ is a given function and the domain $\Omega$ is defined by 
\begin{equation}
\label{Omega}
  \Omega:= \big\{ L^2 < 1 \,;\, L^2= \frac{ |\nabla \varphi|^2 }{c^2 } \big\}, 
\end{equation}
which is an open set in $\R^2$ for $\vp$ sufficiently regular. 

\begin{remark}
\label{Bounded_Omega}
1. It is not difficult to show that, the equation \eqref{edpp} is rotationally invariant.

\medskip
2. In physically realistic problems, one usually imposes global conditions such as:
an incident supersonic flow, or a pseudo-mach number bigger than one at infinity.
Indeed, a typical assumption is
\[
|\nabla \varphi(\xi)| \;\longrightarrow\; U_\infty, 
\quad \text{as } |\xi| \to \infty, 
\]
and the flow is pseudo-supersonic at infinity, that is to say,  
\[
\lim_{|\xi|\to\infty} L^2(\xi) \;>\; 1.
\]
Consequently, the subsonic domain $\Omega$ cannot extend to infinity, which means that, it is bounded. 
\end{remark}

%%%%%%%%%%%%%%%%%
\subsection{Ellipticity Principle}
%%%%%%%%%%%%%%%%%

The main focus of this section is the Ellipticity Principle, which
states that, the interior of the parabolic-elliptic region, $\{L \leq 1\}$,  
of a sufficiently smooth solution of \eqref{edpp} must be elliptic.
This result was first established by Elling, Liu \cite{EllingLiu},
for $\gamma > -1$. Then Chen, Feldman \cite{ChenFeldman}, Chapter 5, extended this result, 
providing a more precise estimate for the case $\gamma \geq 1$.
In absence of such a result, a possibly 
infinite number of parabolic-hyperbolic ``bubbles'' could arise
in the interior of the flow domain. These “bubbles” correspond to regions where the equation \eqref{edpp} changes type, 
switching from elliptic to hyperbolic, with the interface itself being parabolic.
From a physical perspective, such configurations would represent highly unstable patterns of alternating subsonic and supersonic regions within the flow, 
which are not expected to occur in realistic compressible fluids. 
See also \cite{Liu1}, Chapter 13.

\medskip
Here, we study the case $\gamma \ge -1$ closely following the proof by Chen-Feldman. 
We note that items $(i)$ and $(ii)$ of Theorem \ref{Ellipticity_principle_theorem}, below, hold for 
$\gamma > -1$, whereas item $(i)$ also remains valid for $\gamma = -1$. Moreover, we recall the following notation
$\varphi_i= \dfrac{\p \varphi}{\p x_i}$, and for $\sigma \in \R$, 
$O(\sigma)$ denotes any expression such that, $|O(\sigma)| \le C |\sigma|$, where $C$ may depend on $\gamma$.

\begin{theorem}[Ellipticity Principle]
\label{Ellipticity_principle_theorem}
Let $\Omega \subset \R^2$ be an open bounded domain. 

\begin{description}
\item [$(i)$] Suppose $\gamma \ge -1$. Let $\varphi \in C^3(\Omega)$ satisfy (\ref{edpp}) with $L \le 1$ and $\rho>0$ in $\Omega$. Then either $L \equiv 0$ in $\Omega$ or $L$ does not attain its maximum in $\Omega$.

\item [$(ii)$] Suppose $\gamma>-1$. For any $D>0$, there exists $C_0>0$ depending only on $(\gamma,D)$ such that, if $\diam (\Omega) \le D$, 
for any $ \delta \ge 0$, $\hat{c} \ge 1$, and $b \in C^2(\Omega)$ with $|Db|+ \hat{c}|D^2b| \le \frac{\delta}{\hat{c}}$, and for any solution 
$\varphi \in C^3(\Omega)$ of (\ref{edpp}) satisfying $L \le 1$, $\rho(|D\varphi|^2, \varphi)>0$, and $c(|D\varphi|^2,\varphi) \le \hat{c}$ in $\Omega$, then either
$$L^2 \le C_0 \, \delta \quad \text{in $\Omega$}$$
or $L^2+b$ does not attain its maximum  in $\Omega$.
\end{description}
\end{theorem}

\begin{remark}
We remark that Theorem \ref{Ellipticity_principle_theorem} is applied twice in this paper. 
To the best of the authors’ knowledge, both applications appear for the 
first time for that purpose: first in the proof of Theorem \ref{theorem4.3}, and 
second in the proof of Theorem \ref{QPTHM}.
\end{remark}

\begin{proof}
1. First we will prove $(ii)$. We suppose that $\gamma>-1$, $\gamma \neq 1$.
Let $D>0$, $\diam (\Omega) \le D$, $\delta>0$ and $\hat{c}\ge 1$. Fix $b \in C^2(\Omega)$ with $|Db|+\hat{c}|D^2b| \le \frac{\delta}{\hat{c}}$. 
Let $\varphi \in C^3(\Omega)$ solution of (\ref{edpp}) with $L \le 1$, $\rho >0$, $c \le \hat{c}$ in $\Omega$. 
Suppose $L^2+b$ attains its maximum at $\hat{P} \in \Omega$. Let $L^2(\hat{P})> C_1 \delta$, for a constant $C_1 \ge 1$ to be chosen. 

We note that, $L^2(\hat{P})> 0$ implies $|\nabla \varphi (\hat{P})|= L(\hat{P}) \, c(\hat{P})> 0$, and since equation \eqref{edpp} 
is rotationally invariant we can assume, without loss of generality, that 
$\varphi_1(\hat{P})> 0$ and $\varphi_2(\hat{P}) = 0$. Therefore, we have
$$
    0< \varphi_1(\hat{P})= |\nabla \varphi(\hat{P})|= L(\hat{P}) \, c(\hat{P}).
$$
From now on, unless otherwise specified, the equations and inequalities are meant to hold in the maximum point $\hat{P}$, and for $\delta \in \R$.
We also note that, 
\begin{equation}\label{deltaM}
\dfrac{\delta}{L}< 1.
\end{equation}
Indeed, by hypothesis $L^2> C_1 \delta$. Therefore, we have 
$$\dfrac{\delta}{L}< \dfrac{L}{C_1} \le 1,$$
where the last inequality holds because $L \le 1$ and $C_1 \ge 1$.
Moreover, since $L^2+b$ attains its maximum at $\hat{P}$, thus the first and second order derivatives conditions implies that
\begin{equation}
\label{maximumpoint1}
(L^2+b)_1=0, \quad 
(L^2+b)_2=0,
\end{equation}
and
\begin{equation}
\label{maximumpoint2}
(L^2+b)_{11}\le 0, \quad
(L^2+b)_{22}\le 0.
\end{equation}

\medskip
2. Now, we substitute 
$$
L^2 (\xi)= \dfrac{|\nabla\varphi (\xi)|^2}{c^2}, \quad
\varphi_1 (\hat{P})= L(\hat{P})c(\hat{P}), \quad
c^2(\xi)=-(\gamma-1)\left( \varphi + \dfrac{1}{2}| \nabla \varphi |^2 \right) + C
$$
into \eqref{maximumpoint1}, \eqref{maximumpoint2}.
Then, calculating the first and second order derivatives of $c^2(\xi)$, where $c^2$ is given by equation \eqref{eqc2}, we have
\begin{equation*}
\begin{split}
&c^2_1(\xi)= - (\gamma -1)(\varphi_1 + \varphi_1 \varphi_{11} + \varphi_2 \varphi_{21}),\\[5pt]
&c^2_{11}(\xi)= - (\gamma -1)(\varphi_{11} + \varphi_1 \varphi_{111} + \varphi^2_{11} + \varphi_2 
\varphi_{211} + \varphi^2_{21}),\\[5pt]
&c^2_2(\xi)=- (\gamma -1)(\varphi_2 + \varphi_1 \varphi_{12} + \varphi_2 \varphi_{22}),\\[5pt]
&c^2_{22}(\xi)=-(\gamma -1)(\varphi_{22}+ \varphi_1 \varphi_{122} + \varphi^2_{12} + \varphi_2 \varphi_{222} + \varphi^2_{22}).
\end{split}
\end{equation*}
At $\hat{P}$ we have $\varphi_2=0$ and $\varphi_1= L \, c$, and hence 
\begin{equation}\label{derivc2first}
c^2_1=-L \, c (\gamma-1)(1+\varphi_{11}),\quad
c^2_2= - L \,  c (\gamma-1)\varphi_{12},
\end{equation}
and
\begin{equation}\label{derivc2second}
\begin{split}
&c^2_{11}=-(\gamma-1)(L \, c \, \varphi_{111}+\varphi_{11}+\varphi^2_{11}+\varphi^2_{12}),
\\[5pt]
&c^2_{22}=-(\gamma-1)(L \, c \, \varphi_{122}+\varphi^2_{12}+\varphi_{22}+\varphi^2_{22}).
\end{split}
\end{equation}
Similarly, for $L^2$ we obtain
$$
(L^2)_k= \Big( \dfrac{|D \varphi|^2}{c^2}\Big)_k= \dfrac{c^2 (2 \varphi_1 \varphi_{1k}+2 \varphi_2 \varphi_{2k})-(\varphi_1^2+\varphi_2^2)c^2_k}{c^4}.
$$
Thus, at $\hat{P}$,
\begin{equation}
\label{m2k}
(L^2)_k= \dfrac{2  L \, \varphi_{1k}}{c} - \dfrac{L^2 \, c^2_k}{c^2}.
\end{equation}
Analogously, we obtain
$$
\begin{aligned}
&(L^2)_{kk}=
 \left(\dfrac{c^2 (2 \varphi_1 \varphi_{1k}+2 \varphi_2 \varphi_{2k})-(\varphi_1^2+\varphi_2^2)c^2_k}{c^4} \right)_k\\
%&=& \left( \dfrac{2 \varphi_1 \varphi_{1k}+2 \varphi_2 \varphi_{2k}}{c^2}-
%\dfrac{(\varphi_1^2+\varphi_2^2)(c^2)_k}{c^4} \right)_k\\
%&=& \dfrac{2c^2(\varphi_1 \varphi_{1kk}+\varphi^2_{1k}+\varphi_2 \varphi_{2kk}+\varphi^2_{2k})}{c^4}
%- \dfrac{2c^2_k(\varphi_1 \varphi_{1k} + \varphi_2 \varphi_{2k})}{c^4}\\
%&& - \dfrac{c^4c^2_{kk}(\varphi^2_1 + \varphi^2_2)}{c^8}
%- \dfrac{2c^4c^2_k(\varphi_1\varphi_{1k}+\varphi_2\varphi_{2k})}{c^8}\\
%&&+ \dfrac{2c^2\left(c^2_k \right)^2(\varphi^2_1 + \varphi^2_2)}{c^8}\\
&= \dfrac{2(L c\varphi_{1kk}+\varphi^2_{1k}+\varphi^2_{2k})}{c^2}
- \dfrac{2 L \varphi_{1k}c^2_k}{c^3} - \dfrac{L^2c^2_{kk}}{c^2}
- \dfrac{2 L \varphi_{1k}c^2_k}{c^3}
+\dfrac{2 L^2\left(c^2_k \right)^2}{c^4}
\end{aligned}
$$
and at $\hat{P}$, we have 
\begin{equation}\label{m2kk}
(L^2)_{kk}=\dfrac{2(L c\varphi_{1kk}+\varphi^2_{1k}+\varphi^2_{2k})}{c^2}
- \dfrac{4 L \varphi_{1k}c^2_k}{c^3}
- \dfrac{L^2c^2_{kk}}{c^2}
+\dfrac{2 L^2\left(c^2_k \right)^2}{c^4}.
\end{equation}

\medskip
3. Now, substituting (\ref{derivc2first}) and (\ref{m2k}) into (\ref{maximumpoint1}), 
$$
   \dfrac{2 L \varphi_{11}}{c}+ \dfrac{L^3}{c}(\gamma-1)(1+ \varphi_{11})+b_1= 0.
$$
Multiplying by $\dfrac{c}{L}$, it follows that 
\begin{equation}\label{vp1}
[2+(\gamma-1) L^2]\varphi_{11} + L^2(\gamma-1) + \dfrac{c \, b_1}{L}= 0.
\end{equation}
Similarly, we have
%$$\dfrac{2M\varphi_{12}}{c}+ \dfrac{M^3}{c}(\gamma-1)\varphi_{12}+b_2=0,$$
%and, thus,
\begin{equation}\label{vp2}
[2+(\gamma-1) L^2]\varphi_{12}+ \dfrac{c \, b_2}{L}= 0.
\end{equation}

\underline{Claim \# 1:} If $\gamma>-1$, $\gamma \neq 1$, then $2+(\gamma-1) L^2 \neq 0$ and there exists $C(\gamma)>0$, such that
\begin{equation}\label{2gammaM}
\left| \dfrac{1}{2+(\gamma-1) L^2} \right| < C(\gamma).
\end{equation}
\underline{Proof of Claim:} In fact, if $\gamma > 1$, we have $\gamma-1>0$, thus $2 + L^2(\gamma-1) > 2$. 
Therefore, we have
$$0< \dfrac{1}{2+(\gamma-1) L^2} < \dfrac{1}{2}.$$
If $-1< \gamma <1$, since $L^2 \le 1$ and $\gamma -1 < 0$, we have
$(\gamma -1) L^2 \ge \gamma -1$. Therefore,
$$ 2+ (\gamma -1) L^2 \ge \gamma + 1 >0$$
and
$$0< \dfrac{1}{2+(\gamma-1) L^2} < \dfrac{1}{\gamma+1},$$
which completes the proof of the claim.

\medskip
It follows from the claim that
%\begin{eqnarray*}
%\left| \dfrac{c b_i}{M[2+(\gamma-1)M^2]} \right| &\le& C(\gamma) \cdot \dfrac{c|Db|}{M}\\
%&\le& C(\gamma) \cdot \dfrac{1}{M} \hat{c} \cdot \dfrac{\delta}{\hat{c}}\\
%&=& C(\gamma) \dfrac{\delta}{M},
%\end{eqnarray*}
%for $i=1,2$.
%That is,
\begin{equation}
\dfrac{c \, b_i}{ L [2+(\gamma-1) L^2]} = O\Big( \dfrac{\delta}{L} \Big), \quad (i=1,2).
\end{equation}
Therefore, equations (\ref{vp1}) and (\ref{vp2}) becomes
\begin{equation}\label{vp11}
\varphi_{11}= \dfrac{(1-\gamma)L^2}{2+(\gamma-1)L^2}+ O \Big( \frac{\delta}{L}\Big)
\end{equation}
and
\begin{equation}\label{vp12}
\varphi_{12} =O \Big( \frac{\delta}{L} \Big).
\end{equation}

\medskip
4. Now, we will use equation (\ref{edpp})
at $\hat{P}$, in order to obtain $\varphi_{22}$.
Indeed, at $\hat{P}$ we have $\varphi_2= 0$. Then,
$c^2 \varphi_{11} + c^2 \varphi_{22} - \varphi^2_1 \varphi_{11} = \varphi^2_1 - 2c^2$. 
Also $\varphi_1= L \, c$ at $\hat{P}$, thus, 
$c^2 \varphi_{11}+ c^2 \varphi_{22} - L^2c^2 \varphi_{11} = L^2c^2 - 2c^2$. 
Consequently, we have 
\begin{flalign}
\label{vp22}
\varphi_{22}&= (L^2-1)\varphi_{11}+L^2-2\nonumber 
= (L^2-1) \Big[ \dfrac{(1-\gamma) L^2}{2+(\gamma -1) L^2}+ O \Big( \frac{\delta}{L} \Big) \Big] + L^2-2\nonumber 
\\[5pt]
&= \dfrac{(1-\gamma) L^2(L^2-1) + (L^2 - 2)(2+(\gamma-1) L^2)}{2+(\gamma-1) L^2} + O \Big( \frac{\delta}{L}\Big)\nonumber \\[5pt]
%&\varphi_{22}=\dfrac{(1-\gamma)M^2(M^2-1-M^2+2)+2M^2-4}{2+(\gamma-1)M^2} + O \left( \frac{\delta}{M}\right)\nonumber \\[5pt]
%&\varphi_{22}=\dfrac{(1-\gamma)M^2+2M^2-4}{2+(\gamma-1)M^2} + O \left( \frac{\delta}{M}\right)\nonumber \\[5pt]
&=\dfrac{(3-\gamma)L^2-4}{2+(\gamma-1)L^2} + O \left( \frac{\delta}{L}\right).
\end{flalign}

\medskip
\underline{Claim \# 2:} There exists $C(\gamma)>0$, such that
\begin{equation}
\label{Eq_Claim2}
|\varphi_{11}|+|\varphi_{22}| \le C(\gamma).
\end{equation}

\underline{Proof of Claim:} From \eqref{vp11} and \eqref{vp22}, it follows that, 
\begin{eqnarray*}
|\varphi_{11}|+|\varphi_{22}| &\le& \Big| \dfrac{(1-\gamma) L^2}{2+(\gamma-1) L^2} \Big| 
+ \Big| \dfrac{(3-\gamma) L^2-4}{2+(\gamma-1) L^2} \Big| +  \underbrace{O \Big( \frac{\delta}{L}\Big)}_{| O \left( \frac{\delta}{L}\right)| \le C(\gamma) \cdot \frac{\delta}{L}}.
\end{eqnarray*}
We have seen previously that for $\gamma>-1$, $\gamma \neq 1$,
$$
\Big| \dfrac{1}{2+(\gamma-1)L^2} \Big| \le C(\gamma), 
$$
and since $\frac{\delta}{L}< 1$, the thesis of the claim follows. 

\medskip
5. Now, we estimate $\varphi_{111}$. Substituting (\ref{derivc2second}) and (\ref{m2kk}), (with $k=1$), into (\ref{maximumpoint2}), 
we obtain 
\begin{equation*}
\begin{split}
&\dfrac{2(L \, c \varphi_{111} + \varphi^2_{11} + \varphi^2_{12})}{c^2}
+ \dfrac{4 L^2 \varphi_{11}(\gamma-1)(1+\varphi_{11})}{c^2}
\\
&+\dfrac{L^2(\gamma-1)(L \, c\varphi_{111} + \varphi_{11} + \varphi^2_{11} + \varphi^2_{12})}{c^2}
+\dfrac{2 L^4(\gamma-1)^2(1+\varphi_{11})^2}{c^2}
+b_{11}\le 0, 
\end{split}
\end{equation*}
and thus 
\begin{equation*}
\begin{split}
&L \, c \, \varphi_{111}[2 + (\gamma-1) L^2] + \varphi^2_{12}[2+(\gamma-1) L^2]
\\
&+2 \varphi_{11}[\varphi_{11} + 2(\gamma-1) L^2(1+\varphi_{11})]
\\
&+(\gamma-1) L^2(1+\varphi_{11})[\varphi_{11}+2(\gamma-1) L^2(1+\varphi_{11})]
+c^2b_{11} \le 0.
\end{split}
\end{equation*}
Therefore, we get 
\begin{equation}
\label{vp111first}
\begin{aligned}
L \, c \, \varphi_{111} \le  
&- \dfrac{[2\varphi_{11}+(\gamma-1)L^2(1+\varphi_{11})][\varphi_{11}+2(\gamma-1)L^2(1+\varphi_{11})]}{2+(\gamma-1)L^2}
\\
&- \varphi^2_{12}- \dfrac{c^2b_{11}}{2+(\gamma-1)L^2}.
\end{aligned}
\end{equation}
Analogously, we estimate $\varphi_{122}$ due to (\ref{derivc2second}) and (\ref{m2kk}) into (\ref{maximumpoint2}), (for $k=2$),
\begin{equation*}
\begin{split}
&\dfrac{2(L \, c \, \varphi_{122} + \varphi^2_{12}+\varphi^2_{22})}{c^2}
+\dfrac{4 L^2\varphi^2_{12}(\gamma-1)}{c^2}
\\
&+\dfrac{L^2(\gamma-1)(L \, c \, \varphi_{122} + \varphi^2_{12} + \varphi_{22}+\varphi^2_{22})}{c^2}
+\dfrac{2 L^4(\gamma-1)^2\varphi^2_{12}}{c^2}
+b_{22}\le 0, 
\end{split}
\end{equation*}
and thus
\begin{equation*}
\begin{split}
&L \, c \, \varphi_{122}(2 + (\gamma-1) L^2) + \varphi^2_{22}(2+(\gamma-1) L^2)+(\gamma-1) L^2\varphi_{22}
\\
&+2\varphi^2_{12}(1+2(\gamma-1) L^2)
+(\gamma-1) L^2 \varphi^2_{12}(1+2(\gamma-1)L^2)
+c^2 \, b_{22} \le 0.
\end{split}
\end{equation*}
Consequently, 
\begin{equation}\label{vp122first}
L \, c \, \varphi_{122} \le 
-\varphi^2_{22}- \dfrac{(\gamma-1) L^2 \varphi_{22}}{2+(\gamma-1) L^2}
- \varphi^2_{12}[1+2(\gamma-1) L^2] - \dfrac{c^2b_{22}}{2+(\gamma-1) L^2}.
\end{equation}
Also, applying Claim \# 1, for $i=1,2$,
and using $\hat{c}|D^2b| \le \frac{\delta}{\hat{c}}$, we have
\begin{equation}\label{bii}
\dfrac{c^2b_{ii}}{2+(\gamma-1) L^2} = O \Big( \dfrac{\delta}{L} \Big) \quad (i=1,2).
\end{equation}

\medskip
\underline{Claim \#3:} In (\ref{vp111first}) the term, $\varphi_{11} + 2(\gamma-1) L^2(1+\varphi_{11})$, satisfies
\begin{equation}\label{term1}
(i) \qquad |\varphi_{11}+2(\gamma-1) L^2(1+\varphi_{11})| \le C(\gamma),
\end{equation}
and the term, $2\varphi_{11} + (\gamma-1) L^2 (1+\varphi_{11})$, satisfies
\begin{equation}\label{term2}
(ii) \qquad 2 \varphi_{11} + (\gamma-1) L^2 (1+\varphi_{11})=O \Big( \dfrac{\delta}{L}\Big).
\end{equation}
\underline{Proof of Claim:}
Since $L^2 \le 1$ and $|\varphi_{11}| \le C(\gamma)$, we have $(i)$.
Also,
\begin{equation*}
\begin{split}
&2\varphi_{11} + (\gamma-1) L^2 (1+\varphi_{11})
\\
&= \dfrac{2(1-\gamma) L^2}{2+(\gamma-1) L^2}
+ (\gamma-1) L^2 \left[ \dfrac{(1-\gamma) L^2+2+(\gamma-1) L^2}{2+(\gamma-1) L^2}\right] + O \Big( \dfrac{\delta}{L}\Big)
\\
&=O \Big( \dfrac{\delta}{L}\Big).
\end{split}
\end{equation*}
Moreover, it follows from (\ref{vp12}) that
\begin{equation}\label{vp12square}
\varphi_{12}^2 = O\Big(\dfrac{\delta}{L}\Big).
\end{equation}

Finally, we substitute (\ref{2gammaM}), (\ref{bii}), (\ref{term1}), (\ref{term2}) and (\ref{vp12square}) into equation (\ref{vp111first}) to conclude that
\begin{equation}
\label{vp111}
 L \, c \,  \varphi_{111} \le  O\Big(\dfrac{\delta}{L}\Big),
\end{equation}
and similarly,  
substituting (\ref{2gammaM}), (\ref{vp12square}) and (\ref{bii}) into (\ref{vp122first}),
it follows that 
\begin{equation}
\label{ct1}
\begin{aligned}
 c \, L \, \varphi_{122} \le& - \varphi_{22}^2 - \dfrac{(\gamma-1) L^2}{2+(\gamma-1) L^2}\varphi_{22} 
- \underbrace{(1+2(\gamma-1) L^2)\varphi_{12}^2}_{O\Big(\frac{\delta}{L}\Big)} - \underbrace{\dfrac{c^2b_{22}}{2+(\gamma-1) L^2}}_{O\Big(\frac{\delta}{L}\Big)} 
\\
 \le & \dfrac{((\gamma-3) L^2+4)(2 L^2-4)}{(2+(\gamma-1) L^2)^2} + O\Big(\dfrac{\delta}{L}\Big).
\end{aligned}
\end{equation}

\medskip
6. To follow, we claim that there exists $\tilde{C}(\gamma)> 0$, such that
\begin{equation}
\label{termvp122}
\dfrac{((\gamma-3) L^2+4)(2 L^2-4)}{(2+(\gamma-1) L^2)^2} \le - \tilde{C}(\gamma).
\end{equation}
Indeed, if $\gamma>1$ then $\gamma+1>2$. Also $L^2 \le 1$ implies that $2(2-L^2)> 2$. Then, we have 
$(\gamma+1)(4-2 L^2) > 4$ and thus
\begin{equation}
\label{a3}
- \dfrac{1}{4} (4-2 L^2)^2 <  \dfrac{-4}{(\gamma+1)^2}.
\end{equation}
Moreover, 
since $\gamma>1$, we have
$(\gamma-3) L^2 + 4 >  4 - 2 L^2$. 
Multiplying this inequality by $-(4-2 L^2)<0$, we conclude that 
\begin{equation}\label{a4}
-((\gamma-3) L^2+4)(4-2 L^2) < - (4-2 L^2)^2.
\end{equation}
We have seen previously that if $\gamma> 1$, then 
\begin{equation}\label{a5}
\dfrac{1}{2+(\gamma-1) L^2}\le \dfrac{1}{2}, 
\end{equation}
and using (\ref{a4}) and (\ref{a5}), we obtain
\begin{equation}
\label{a6}
\dfrac{((\gamma-3) L^2+4)(2 L^2-4)}{(2+(\gamma-1) L^2)^2} 
\le - \dfrac{1}{4}(4-2 L^2)^2 \underbrace{\le}_{\mbox{by}\,(\ref{a3})} -\dfrac{4}{(\gamma+1)^2}.
\end{equation}

Now, for $-1 < \gamma < 1$ we have from the proof of Claim \# 1,
\begin{equation}
\label{b1}
0 < \dfrac{1}{(2+(\gamma-1) L^2)^2} \le \dfrac{1}{(\gamma+1)^2}.
\end{equation}
Since $0 \le L^2 \le 1$, we obtain
\begin{equation}\label{b2}
(4-2 L^2) \ge 2.
\end{equation}
Multiplying $L^2 \le 1$ by $\gamma-3< 0$ (since $\gamma < 1$), and adding $4$ we have
\begin{equation}\label{b3}
(4+ L^2(\gamma-3)) \ge \gamma+1>0. 
\end{equation}
From (\ref{b2}) and (\ref{b3}), 
\begin{equation}\label{b4}
(4+L^2(\gamma-3)) (2L^2-4) \le - 2 (\gamma+1), 
\end{equation}
and from (\ref{b4}) and (\ref{b1}),
\begin{equation}
\label{b5}
 \dfrac{(4+ L^2(\gamma-3))(2 L^2-4)}{(2+(\gamma-1) L^2)^2} \le - \dfrac{2(\gamma+1)}{(\gamma+1)^2}= -\dfrac{2}{\gamma+1}.
\end{equation}
Moreover, it follows from (\ref{a6}) and (\ref{b5}) that, for $\gamma>-1$,
$$
    \dfrac{((\gamma-3) L^2+4)(2 L^2-4)}{(2+(\gamma-1) L^2)^2} \le - \underline{C}(\gamma).
$$
Lastly, by combining (\ref{ct1}) and (\ref{termvp122}), we deduce
\begin{equation}
\label{5.2.7}
  c \, L \, \varphi_{122} \le - \tilde{C}(\gamma) + O\Big( \dfrac{\delta}{L} \Big), 
\end{equation}
where this inequality holds for $\gamma>-1$, ($\gamma\neq 0,1$).

\medskip
7. Finally, to reach a contradiction, we establish the following estimate
$$
    c \, L \, \varphi_{122} \ge O\Big( \dfrac{\delta}{L} \Big). 
$$
To begin, let us differentiate equation (\ref{edpp}) with respect to $\xi_1$ to obtain 
\begin{flalign*}
%& (c^2(\varphi_{11}+\varphi_{22})- \varphi_1^2 \varphi_{11} - 2 \varphi_1 \varphi_2 \varphi_{12} - \varphi_2^2 \varphi_{22})_1 = (\varphi_1^2 + \varphi_2^2 - 2c^2)_1 \\[10pt]
& c^2(\varphi_{111}+ \varphi_{122})+(\varphi_{11}+ \varphi_{22})\overbrace{(-(\gamma-1)(\varphi_1+ \varphi_1 \varphi_{11}))}^{=(c^2)_1(\hat{P})}
\\[5pt]
& - \varphi_1^2 \varphi_{111}- \varphi_{11} 2 \varphi_1 \varphi_{11} - 2 \varphi_1 \varphi_2 \varphi_{112} - 2 \varphi_{12}(\varphi_1 \varphi_{12}+ \varphi_2 \varphi_{11})- \varphi_2^2 \varphi_{221} 
\\[5pt]
& - 2 \varphi_2 \varphi_{21} \varphi_{22} =  2 \varphi_1 \varphi_{11} + 2 \varphi_2 \varphi_{12} - 2 \underbrace{(-(\gamma-1)(\varphi_1+ \varphi_1 \varphi_{11}))}_{(c^2)_1(\hat{P})}.
\end{flalign*}
Then, using $\varphi_1= L \, c$ and $\varphi_2= 0$, it follows that
\begin{flalign*}
& (c^2- L^2 c^2) \varphi_{111} + c^2 \,  \varphi_{122}= (\gamma-1) L \, c(1+\varphi_{11}) \Delta \varphi + 2 L \, c \, \varphi_{11}^2  
\\[5pt]
&+ 2 L \, c  \, \varphi_{11} + 2 L \, c \, (\gamma-1) + 2 L \, c \, \gamma \, \varphi_{11} - 2 L \, c \, \varphi_{11} + 2 L \, c \, \varphi_{12}^2.
\end{flalign*}
%Recall that
Using $\varphi^2_{12} = O \Big( \frac{\delta}{L}\Big)$
and multiplying by $\frac{L}{c}$, we have
\begin{eqnarray}
\label{vp122contrad}
%&& (c^2-M^2c^2) \varphi_{111} + c^2 \varphi_{122}= \nonumber \\[5pt]
%&& (\gamma-1) Mc(1+ \varphi_{11}) \Delta \varphi + 2Mc \varphi_{11}^2+2Mc(\gamma-1) + 2Mc \gamma \varphi_{11} + McO \left( \dfrac{\delta}{M} \right)\nonumber \\[15pt]
%&& (c-M^2c) \varphi_{111} + c\varphi_{122}=\nonumber\\[5pt]
%&& (\gamma-1) M(1+ \varphi_{11}) \Delta \varphi + 2M \varphi_{11}^2+2M(\gamma-1) + 2M\gamma \varphi_{11} + MO \left( \dfrac{\delta}{M} \right)\nonumber\\[15pt]
&& c \, L \, \varphi_{122}=(L^2-1) L \, c \, \varphi_{111} + (\gamma-1) L^2(1+\varphi_{11}) \Delta \varphi + \nonumber
\\[5pt]
&& \qquad  \qquad + 2 L^2 \, \varphi_{11}^2 + 2 L^2 \, (\gamma-1) + 2 L^2 \, \gamma \, \varphi_{11} + L^2  O\Big( \dfrac{\delta}{L} \Big).
\end{eqnarray}
It follows from (\ref{vp111}) that, 
$(L^2-1) L \, c \, \varphi_{111} \ge  O\Big( \dfrac{\delta}{L} \Big)$.
Therefore, we obtain 
$$
   c \, L \, \varphi_{122} \ge L^2 ((\gamma-1) (1+ \varphi_{11} ) \Delta \varphi + 2 \varphi_{11}^2 + 2 \gamma \varphi_{11} + 2 (\gamma-1)) + O\Big( \dfrac{\delta}{L} \Big).
$$
Substituting (\ref{vp11}) and (\ref{vp22}) in the above equation, 
it follows that
\begin{flalign*}
c \, L \, \varphi_{122}  \ge & L^2 \Big\{ (\gamma-1) \Big[ \dfrac{2+(\gamma-1) L^2+(1-\gamma) L^2}{2+(\gamma-1) L^2} \Big]  
\\[5pt]
&  \times \Big[ \dfrac{(1-\gamma) L^2}{2+(\gamma-1) L^2} + \dfrac{[(1-\gamma)+2] L^2-4}{2+(\gamma-1) L^2} \Big]  
+ \dfrac{2(1-\gamma)^2 L^4}{(2+(\gamma-1) L^2)^2}  
\\[5pt] 
& + \dfrac{2\gamma(1-\gamma) L^2}{2+(\gamma-1) L^2} + 2(\gamma-1)\Big\}+ O\Big( \dfrac{\delta}{L} \Big).
\end{flalign*}
The expression between braces in the inequality above cancels out. 
Thefore, we have
\begin{equation}
\label{5.2.8}
c \, L \, \varphi_{122} \ge O\Big( \dfrac{\delta}{L} \Big).
\end{equation}

\medskip
8. Equations (\ref{5.2.7}) and (\ref{5.2.8}) lead to a contradiction provided that $\frac{\delta}{L}$ is sufficiently small, that is, 
for sufficiently large $C_1 \geq 1$ depending only on $\gamma$.
Therefore, if $L^2 + b$ attains its maximum at $\hat{P} \in \Omega$, there exists $C_1 \in \mathbb{R}$, such that, 
$L^2(\hat{P}) \le C_1 \delta$. 
Consequently, if $\xi \in \Omega$, then
\begin{eqnarray*}
%M^2(\xi)+b(\xi)&\le & M^2(\hat{P})+b(\hat{P})\\
L^2(\xi) & \le & L^2(\hat{P}) + b(\hat{P}) - b(\xi)
\\
&\le & L^2(\hat{P})+ \mathrm{osc}_{\Omega} (b)
\\
& \le & \big(C_1 + \diam (\Omega) \big) \, \delta 
=: C_0 \, \delta,
\end{eqnarray*}
where we have used: $\mathrm{osc}_{\Omega}(b) \le \dfrac{\delta}{\hat{c}} \diam (\Omega) \underbrace{\le}_{\hat{c}\ge 1} \delta \diam (\Omega)$.

\medskip
9. Now, let us consider the case $\gamma=1$. We recall that, the pressure law is given by
$$p(\rho)=a^2(\rho - \underline{\rho}), \qquad \rho \ge \underline{\rho}>0,$$
and the sound speed is
$c^2=p'(\rho)=a^2$. 
So, we have $c= a$ and
$$
(c^2)_k= 0, \quad
(c^2)_{kk}= 0, \quad 
(k=1, 2).
$$
Therefore, equations (\ref{m2k}) and (\ref{m2kk}) become
\begin{equation}
\label{mkg1}
(L^2)_k= \dfrac{2 L \, \varphi_{1k}}{a}
\end{equation}
and
\begin{equation}
\label{mkkg1}
(L^2)_{kk}= \dfrac{2(L \, a \, \varphi_{1kk} + \varphi^2_{1k} + \varphi^2_{2k})}{a^2}.
\end{equation}
Due to $(L^2+b)_k= 0$, and (\ref{mkg1}) it follows that
\begin{equation}\label{vp1kg1}
\varphi_{1k}= O\Big( \dfrac{\delta}{L} \Big), \quad (k=1, 2).
\end{equation}
Using $c= a$, $\varphi_1= L \, a$ and $\varphi_2= 0$ at (\ref{edpp}), we have
$a^2 \varphi_{11} + a^2 \varphi_{22} - L^2 a^2 \varphi_{11} = L^2 a^2 - 2 a^2$, and thus 
\begin{equation}
\label{vp22g1} 
\varphi_{22}= L^2 - 2 + O\Big( \dfrac{\delta}{L} \Big).
\end{equation}
Substituting (\ref{mkkg1}) into $(L^2+b)_{11} \le 0$ and $(L^2+b)_{22} \le 0$, respectively, it follows that
\begin{equation}
\label{vp111g1}
 \dfrac{2(L \, a \, \varphi_{111} + \varphi^2_{11} + \varphi^2_{12})}{a^2} + b_{11} \le 0, \quad \text{thus} \quad
 \varphi_{111} \le O\Big( \dfrac{\delta}{L}\Big)
\end{equation}
also
\begin{equation}
\label{vp122g1.1}
 \dfrac{2(L \, a\, \varphi_{122} + \varphi^2_{12}+ \varphi^2_{22})}{a^2} + b_{22} \le 0, \quad \text{so} \quad
%&& Ma \varphi_{122} + (M^2-2)^2 + O \left( \dfrac{\delta}{M}\right) \le 0 \nonumber \\
\varphi_{122} \le - \dfrac{(L^2-2)^2}{L \, a} + O \Big( \dfrac{\delta}{L} \Big).
\end{equation}
Now, substituting $c=a$ into (\ref{edpp}) and diﬀerentiating it with respect to $\xi_1$, we obtain
\begin{eqnarray*}
&& a^2(\varphi_{111}+\varphi_{122}) - \varphi^2_1 \varphi_{111} - 2 \varphi_1 \varphi^2_{11}
-2 \varphi_1 \varphi_2 \varphi_{112} - 2 ( \varphi_1 \varphi_{12} + \varphi_2 \varphi_{11}) \varphi_{12}\\
&& - \varphi^2_2 \varphi_{122} -2 \varphi_2 \varphi_{12} \varphi_{22} = 2 \varphi_1 \varphi_{11} + 2 \varphi_2 \varphi_{12}.
\end{eqnarray*}
Using $\varphi_2=0$, we have
\begin{equation}
\label{vp122g1.2}
%&& a(\varphi_{111} + \varphi_{122}) - M^2 a \varphi_{111} - 2 M  \varphi^2_{11} - 2M \varphi^2_{12} = 2M \varphi_{11} \nonumber \\
%&& a(1-M^2) \varphi_{111} + a \varphi_{122} = O \left( \dfrac{\delta}{M}\right) \nonumber \\
 a \varphi_{122} = - a (1-L^2) \varphi_{111} + O \Big( \dfrac{\delta}{L} \Big), \quad \text{then} \quad 
 \varphi_{122} \ge O \Big( \dfrac{\delta}{L}\Big),
\end{equation}
from which we have a similar contradiction as before. Then, the proof of item $(ii)$ is finished.  

\medskip
10. To conclude, we prove $(i)$. Suppose that $L$ is not identically zero in $\Omega$ and attains its maximum in $\hat{P} \in \Omega$. Therefore, we have 
\begin{equation}
L^2> 0 \quad \mbox{at} \quad \hat{P}.
\end{equation}
As before, we have $|\nabla \varphi|= L \, c> 0$ at $\hat{P}$, and we suppose $\varphi_1> 0$ and $\varphi_2= 0$.
Now, we have
\begin{equation}
\label{mderiv}
(L^2)_k= 0, \quad 
(L^2)_{kk}> 0,
\quad (k=1, 2).
\end{equation}
If $\gamma > -1$ and $\gamma \neq 1$, then we apply equations (\ref{m2k}), (\ref{m2kk}) and proceed as in the proof of $(ii)$.
We see that, all the equations obtained previously remain valid just by replacing the term $O\left(\frac{\delta}{L}\right)$ by $0$.
Therefore, there exists $\tilde{C}(\gamma)>0$, such that, $c \, L \, \varphi_{122} \le - \tilde{C}(\gamma)$, and also we have 
$c \, L \, \varphi_{122} \ge 0$, which is a contradiction.

\medskip
Now, if $\gamma= 1$ by (\ref{vp1kg1}) and (\ref{vp22g1}), we have
$$
\varphi_{11}=\varphi_{12}= 0, 
\quad 
\varphi_{22}= L^2-2.
$$
And by (\ref{vp111g1}), (\ref{vp122g1.1}) and (\ref{vp122g1.2}), it follows that, $\varphi_{111} \le 0$ and 
we obtain a contradiction
$$
 \varphi_{122} \le - \dfrac{(L^2-2)^2}{L \, a}, 
\quad \varphi_{122} \ge 0.
$$

\medskip
Now, for $\gamma=-1$ we proceed analogously to the previous case. 
Then, by equations (\ref{mderiv}), (\ref{m2k}) and (\ref{m2kk}), we have
\begin{equation}
\label{eq1gmenos1}
2(1- L^2) \varphi_{11} - 2 L^2= 0,
\end{equation}
and
\begin{equation}
\label{eq2gmenos1}
2(1- L^2) \varphi_{12}= 0.
\end{equation}
It follows by equation (\ref{eq1gmenos1}) that $L^2<1$ at $\hat{P}$.
So, isolating the terms $\varphi_{11}$ and $\varphi_{12}$ in equations (\ref{eq1gmenos1}) and (\ref{eq2gmenos1}), we deduce
$$
\varphi_{11}= \dfrac{L^2}{1 - L^2},
\qquad 
\varphi_{12}= 0.
$$
Replacing $\gamma$ by $-1$ and $O\Big( \dfrac{\delta}{L}\Big)$ by $0$ at (\ref{vp22}), we have
$$
\varphi_{22}= \dfrac{4 L^2-4}{2-2 L^2}= -2. 
$$
Now, from (\ref{vp111first}) it follows that
$$
L \, c \, \varphi_{111} \le - \dfrac{[2 \varphi_{11} - 2 L^2(1 + \varphi_{11})][\varphi_{11} - 4 L^2(1+ \varphi_{11})]}{2(1 - L^2)}, 
$$
and since
$$2 \varphi_{11} - 2 L^2 (1+ \varphi_{11}) = \dfrac{2 L^2}{1 - L^2} - 2 L^2 \left( 1 + \dfrac{L^2}{1 - L^2} \right)= 0,$$
we obtain 
\begin{equation}
\label{eq4gmenos1}
L \, c \, \varphi_{111} \le 0.
\end{equation}
From (\ref{vp122first}) we have
\begin{equation}
\label{eq4gmenos1new}
L \, c \,  \varphi_{122} \le  - \varphi^2_{22} + \dfrac{2 L^2 \varphi_{22}}{2(1 - L^2)}
 \le  -4 - \underbrace{\dfrac{2 L^2}{1 - L^2}}_{\le 0} 
 \le  -4. 
\end{equation}
On the other hand, from (\ref{vp122contrad})
\begin{equation}
\label{eq5gmenos1}
\begin{aligned}
L \, c \,  \varphi_{122} &= (L^2-1) L \, c \, \varphi_{111} - 2 L^2 (1+ \varphi_{11})(\varphi_{11} + \varphi_{22})
\\[5pt]
&\quad + 2 L^2 \varphi^2_{11} - 4 L^2 - 2 L^2 \varphi_{11} 
\\[5pt]
&= \underbrace{(L^2-1) L \, c \varphi_{111}}_{\ge 0, \, \mbox{by} \, (\ref{eq4gmenos1})} - 2 L^2 \Big( 1+ \dfrac{L^2}{1 - L^2} \Big) \Big( \dfrac{L^2}{1 - L^2} -2 \Big)
\\[5pt]
&\quad + 2 L^2 \dfrac{L^4}{(1-L^2)^2} - 4 L^2 - \dfrac{2 L^4}{1-L^2} 
\\[5pt]
& \ge  \dfrac{-2 L^2(3 L^2-2)}{(1- L^2)^2} + \dfrac{2 L^6}{(1- L^2)^2}
\\[5pt]
&\quad \dfrac{-4 L^2(1-2 L^2 + L^4)}{(1- L^2)^2} - \dfrac{2 L^4(1- L^2)}{(1 - L^2)^2}= 0.
%\\[5pt]
%& \ge  \big\{\!-6 L^2 + 4 + 2 L^4 - 4  + 8 L^2 - 4 L^4 - 2 L^2 + 2 L^4 \big\} \dfrac{L^2}{(1 - L^2)^2}= 0.
\end{aligned}
\end{equation}
Finally, equations (\ref{eq4gmenos1new}) and (\ref{eq5gmenos1}) 
gives us a contradiction, which completes the proof.
\end{proof}

%%%%%%%%%%%%%%%%%%%%%
\subsection{Existence of potential flows} 
%%%%%%%%%%%%%%%%%%%%%

The aim of this section is to show that the Dirichlet problem
\begin{equation}\label{dirichlet}
\left\{
\begin{array}{l}
c^2 \Delta \varphi - D^2  \varphi  \nabla \varphi \cdot \nabla \varphi = | \nabla \varphi|^2 - 2c^2 \quad \text{in $\Omega$}, 
\\[5pt]
\varphi = \varphi_b \quad \text{on $\partial \Omega$},
\end{array}
\right.
\end{equation}
has a solution $\varphi \in C^3 (\Omega) \cap C^{2,\alpha} (\overline{\Omega})$, where $c^2$ is given by equation \eqref{eqc2}, 
and 
%$\partial \Omega \in C^{3, \alpha}$, 
$\varphi_b \in C^{2, \alpha} (\overline{\Omega})$, $\alpha \in (0,1)$. 
In fact, since the boundary condition is given a priori, for instance, by some potential incident flux or geometrical condition, 
we are going to assume that, there exists an extension $\tilde{\varphi}_b \in C^{2,\alpha}(\R^2)$ such that, 
$\varphi_b= \tilde{\varphi}_b$ in $\overline{\Omega}$. 
This assumption further simplifies the construction of a solution to this challenging free-boundary problem, since the domain 
$\Omega$ given by \eqref{Omega} depends on the solution $\vp$ itself. Moreover, we assume that $\Omega$ is bounded, 
see Remark \ref{Bounded_Omega} item 2.  

\medskip
To follow, we observe that substituting (\ref{eqc2}) into (\ref{dirichlet}), $\varphi$ 
should satisfies the following system
\begin{equation}
\label{qigual0}
\left\{
\begin{array}{l}
Q \varphi = 0 \quad \mbox{in} \quad \Omega,
\\
\varphi = \tilde{\varphi}_b \quad \mbox{on} \quad \partial \Omega,
\end{array}
\right.
\end{equation}
where $Q$ is the elliptic operator given by
\begin{equation}\label{qelliptic}
Q \varphi = (c^2 - \varphi_1^2) \varphi_{11} - 2 \varphi_1 \varphi_2 \varphi_{12} + (c^2 - \varphi_2^2)\varphi_{22} - \gamma \nabla \varphi \cdot \nabla \varphi - 2 (\gamma-1) \varphi.
\end{equation}
The operator $Q$ can be rewritten as
$$
%\begin{equation}\label{qgeral}
Q \varphi = A(\gamma, \varphi, \nabla \varphi) :  D^2 \varphi + B(\gamma, \nabla \varphi) \cdot \nabla \varphi + C(\gamma, \varphi) \varphi,
%\end{equation}
$$
where the matrix 
$$
%\begin{equation}\label{matrix_a}
A = [ A_{ij}] = \left[ \begin{array}{cc}
c^2 - \varphi_1^2 & - \varphi_1 \varphi_2 \\
- \varphi_1 \varphi_2 & c^2 - \varphi_2^2
\end{array} \right ]
%\end{equation}
$$
has the eigenvalues
\begin{equation}\label{eigenmatrix_a}
0 \le \lambda = c^2 - |\nabla \varphi|^2 < \Lambda = c^2.
\end{equation}
Since $\lambda$ may be $0$ at $\partial \Omega \equiv \{L^2= 1\} \equiv \{c^2 =| \nabla \varphi|^2 \}$, 
namely sonic-line, 
the operator $Q$ is not uniformly elliptic. 
In order to consider a uniformly elliptic operator we regularize the equation $\eqref{qelliptic}$, that is, we consider for each $\varepsilon> 0$,
\begin{equation}\label{qregularized}
\left\{
\begin{array}{l}
Q_{\varepsilon} \varphi_\varepsilon = 0 \quad \mbox{in} \quad \Omega, 
\\
\varphi_\varepsilon = \tilde{\varphi}_b \quad \mbox{on} \quad \partial \Omega,
\end{array}
\right.
\end{equation}
where $\Omega= \Omega(\vp_\ve)$ is similarly defined as \eqref{Omega} and $Q_{\varepsilon} \varphi_\ve= Q \varphi_\ve+\varepsilon \Delta \varphi_\ve$, 
that is to say,
\begin{equation}\label{qegeral}
Q_{\varepsilon} \varphi_\ve = A_{\varepsilon}(\gamma, \varphi_\ve, \varphi_{\ve 1}, \varphi_{\ve 2}) :  D^2\varphi_\ve + B(\gamma, \nabla \varphi_\ve) \cdot \nabla \varphi_\ve + C(\gamma, \varphi_\ve) \varphi_\ve,
\end{equation}
where
$$
%\begin{equation}\label{matrix_ae}
A_{\varepsilon} = [(A_{\varepsilon})_{ij}] = A + \varepsilon I = \left[ \begin{array}{cc}
c^2 - \varphi_{\ve 1}^2 + \varepsilon & - \varphi_{\ve 1} \varphi_{\ve 2}
\\
- \varphi_{\ve 1} \varphi_{\ve 2} & c^2 - \varphi_{\ve 2}^2 + \varepsilon
\end{array} \right].
%\end{equation}
$$
Therefore, the eigenvalues of $A_{\varepsilon}$ are
\begin{equation}\label{eigenmatrix_avarepsilon}
0 < \lambda = c^2 - |\nabla \varphi_\ve|^2 + \varepsilon < \Lambda = c^2 + \varepsilon,
\end{equation}
and, thus, $Q_{\varepsilon}$ is uniformly elliptic. 

\medskip
The strategy is to show that, for any $\ve> 0$ fixed, the system \eqref{qregularized} admits a solution $\vp_\ve$, which should 
be uniformly bounded with respect to $\ve$. Then, under uniform estimates on the gradient of $\vp_\ve$, we may pass to the limit as $\ve \to 0$, (where
the Ellipticity Principle plays a crucial role), thereby establishing the existence of a solution to system \eqref{qigual0}. 

\begin{theorem}
\label{theorem4.1}
Given $\tilde{\varphi}_b \in C^{2,\alpha}$ and $\gamma\in [-1, + \infty) \backslash  \{0\}$, 
the Dirichlet problem \eqref{qregularized} has a solution $\varphi_\varepsilon \in C^{2,\alpha}(\overline{\Omega})$.
\end{theorem}

\begin{proof}
1. For simplicity of notation, we drop the subscript $\varepsilon$ in $\varphi_\varepsilon$, and simply write $\varphi$. 
We proceed to apply the Schauder Fixed Point Theorem. 
For any $\tilde{w} \in C^{2,1}(\R^2)$ we consider  
$$
  \Omega= \big\{\xi \in \R^2  \,;\, |\nabla \tilde{w}(\xi)|^2 < c^2(\tilde{w}(\xi)) \big\},  
$$
and let $M> 0$ be a constant independently of $\ve>0$. 
Then, we define 
\begin{equation}
\label{Set_S}
S:= \Big\{ w= \tilde{w}|_{\overline{\Omega}} \, / \,  w|_{\partial \Omega} = \tilde{\varphi}_b ,\, |w| \le M  \Big\}. 
\end{equation}
For any $w \in S$, fixed, we define the linear operator $L_{\varepsilon}$ obtained from $Q_{\varepsilon}$ given by
$$
\begin{aligned}
L_{\varepsilon} \varphi&= (c^2 - w_1^2 + \varepsilon) \varphi_{11} - 2 w_1 w_2 \varphi_{12} + (c^2 - w^2_2 + \varepsilon) \varphi_{22}
\\[5pt]
& - \gamma \nabla w \cdot \nabla \varphi - 2 (\gamma -1) \varphi,
\end{aligned} 
$$
and we consider the Dirichlet problem 
\begin{equation}\label{Lepdirichlet}
\left\{
\begin{aligned}
    L_{\varepsilon} \varphi &= 0  \quad \,\,\, \mbox{in} \quad \Omega, 
    \\
     \varphi &= \tilde{\varphi}_b \quad \mbox{on} \quad \partial \Omega. 
\end{aligned}
\right.
\end{equation}
Then, $L_{\varepsilon}$ is a uniformly elliptic linear operator defined in $\Omega= \Omega(w)$ (a known bounded open set),   
whose coefficients lies in $C^{1,\alpha}(\overline{\Omega})$, sufficiently consider $C^{0,\alpha}$, 
and consequently \eqref{Lepdirichlet} is a standard Dirichlet value problem. 

\medskip
2. Now, if $\gamma \geq 1$, then the zero order coefficient in $L_\ve$, that is, $-2(\gamma-1)$ is non-positive. 
Therefore, we can apply Theorem 6.14 in \cite{Gilbarg} to conclude that, the Dirichlet problem \eqref{Lepdirichlet} has a unique solution $\varphi \in C^{2, \alpha}(\overline{\Omega})$,
and thus the operator $w \mapsto T(w)= \vp$ is well defined,
where $\varphi$ is the unique solution of the Dirichlet problem (\ref{Lepdirichlet}). We notice that, since $w \in S$, thus
$$\varphi=T(w) \in C^{2,\alpha} \left( \overline{\Omega}\right)  \hookrightarrow C^{2,1} \left( \overline{\Omega} \right).$$ 
Also, $\varphi|_{\partial \Omega}= \tilde{\varphi}_b$, and it follows from the maximum principle, see Corollary 3.2 in \cite{Gilbarg}, 
that there exists a positive constant $M$, such that $|\varphi| \le M$. Consequently, $T(S) \subset S$.

\medskip
On the other hand, if $\gamma< 1$, $(\gamma \neq 0)$, then the zero order coefficient in $L_\ve$ 
becomes positive, and the previously established result no longer follows directly, see p. 107 in \cite{Gilbarg}.
First, we use the Fredholm Alternative, (see Theorem 6.15 in \cite{Gilbarg}), to prove that the 
Dirichlet problem \eqref{Lepdirichlet} has a unique solution. Since the linear operator $L_\ve$ is uniformly 
elliptic, it is sufficient to prove its injectiveness to apply the Fredholm aternative.
Define $L^{\varepsilon}_-:= - L_{\varepsilon}$, that is, 
$$
\begin{aligned}
L^{\varepsilon}_- \varphi&= -A^\ve(\xi) :  D^2\varphi + \gamma \nabla w \cdot \nabla \varphi + 2 (\gamma -1) \varphi,
\end{aligned} 
$$
where
$$
A^{\varepsilon} =  
\left[ \begin{array}{cc}
(c^2 - w_{1}^2 + \varepsilon) &  -w_{1} w_{2}
\\
 -w_{1} w_{2} & (c^2 - w_{2}^2 + \varepsilon)
\end{array} \right], 
$$
and equivalently we show that $L^{\varepsilon}_-$ is injective. Conveniently, we write 
$L^{\varepsilon}_-$ in divergence form, that is to say,
\begin{equation}
\label{LEPSILONM}
L^{\varepsilon}_- \varphi = - \partial_j \big(a^\ve_{ij}(\xi) \partial_i \varphi \big) + b_i(\xi) \partial_i \varphi+ 2(\gamma-1)\varphi,
\end{equation}
where for all $\xi \in \Omega$, 
$$
\begin{aligned}
    b&= {\rm div} A^\ve + \gamma \, \nabla w
\\[5pt]
&= -(\gamma - 1) \, \nabla w - \gamma \, D^2 w \nabla w - (\Delta w) \nabla w + \gamma \, \nabla w
\\[5pt]
&= \nabla w - \gamma \, D^2 w \nabla w - (\Delta w) \nabla w. 
\end{aligned}    
$$     
To follow, we consider, for some $\kappa_0> 0$ to be chosen a posteriori,
the bilinear form $B[u,v]$  associated with the shifted operator
$L^{\varepsilon}_- + \kappa_0 I$, 
defined for each $u,v \in H_0^1(\Omega)$ by
$$
B[u,v]:= \int_{\Omega}a_{ij}(\xi) \partial_i u \, \partial_j v \, d\xi 
+ \int_{\Omega} b_i(\xi) \partial_i u \, v \, d\xi + (2(\gamma-1) + \kappa_0) \int_{\Omega}  u \, v  \, d\xi. 
$$
The continuity of $B$ follows easily from the regularity of their coefficients, 
and also the boundedness on $H_0^1(\Omega)$.
We claim that it is coercive in $H_0^1(\Omega)$.
Indeed, for $\vp \in H_0^1(\Omega)$, we have 
$$
\begin{aligned}
B[\vp,\vp]& \geq \int_{\Omega} (c^2 - |\nabla w|^2 + \varepsilon) \, |\nabla \vp|^2 \, d\xi 
+ \int_{\Omega} \vp \, b \cdot  \nabla \vp \, d\xi 
\\[5pt]
& + (2(\gamma-1) + \kappa_0) \int_{\Omega} \vp^2  \, d\xi
\\[5pt]
&\geq \ve \, \int_\Omega | \nabla \vp|^2 \, d\xi 
+ \int_{\Omega} \vp \, \nabla w \cdot  \nabla \vp \, d\xi 
- \gamma \int_{\Omega} \vp \, D^2 w \nabla w \cdot  \nabla \vp \, d\xi
\\[5pt]
&- \int_{\Omega} \vp \, \Delta w \nabla w \cdot  \nabla \vp \, d\xi
+( 2(\gamma -1) + \kappa_0) \int_{\Omega}  \vp^2  \, d\xi. 
\end{aligned}
$$
Let us consider the mixed terms, that is, 
$$
I= \gamma \int_{\Omega} \vp \, D^2 w \nabla w \cdot  \nabla \vp \, d\xi
+  \int_{\Omega} \vp \, \Delta w \nabla w \cdot  \nabla \vp \, d\xi
-  \int_{\Omega} \vp \, \nabla w \cdot  \nabla \vp \, d\xi, 
$$
and define 
$$
\tilde{M}= \|\nabla w\|_{L^\infty} + \|D^2 w\|_{L^\infty} + \|\Delta w\|_{L^\infty}. 
$$
 Then, there exists $C= C(\gamma, \tilde{M})> 0$, such that, 
 $$
 \begin{aligned}
 | I | &\leq C \, \int_\Omega |\vp| \, |\nabla \vp| \, d\xi
 \\[5pt]
 &\leq C \, \mu_0 \int_\Omega  |\nabla \vp|^2 \, d\xi + \frac{C}{\mu_0} \int_\Omega  |\vp|^2 \, d\xi, 
 \end{aligned} 
 $$
 for some positive constant $\mu_0$, with $C \mu_0= \ve / 2$.  
 Therefore, we obtain 
 \begin{equation}
 \label{eq:coercive-est}
\begin{aligned}
B[\vp,\vp]
&\geq \ve \, \int_\Omega | \nabla \vp|^2 \, d\xi 
- \frac{\ve}{2} \int_{\Omega} |\nabla \vp|^2 \, d\xi 
\\[5pt]
&+\big ( 2(\gamma -1) - \frac{2 C^2}{\ve} + \kappa_0 \big) \int_{\Omega}  |\vp|^2  \, d\xi
 \geq \frac{\ve}{2} \, \|\vp\|^2_{H^1_0(\Omega)}, 
\end{aligned}
\end{equation}
where we have set
$$ 
   k_0=  - 2 (\gamma - 1)+ \frac{2 C^2}{\ve} + \frac{\ve}{2}. 
$$ 

\medskip
3. Injectivity of $L_-^\varepsilon$ via the shifted operator.
Let $T_{\kappa_0}:H_0^1(\Omega)\to H^{-1}(\Omega)$ 
be the operator associated with $B$, that is to say, 
\[
\langle T_{\kappa_0}u, v\rangle_{H^{-1},H_0^1}=B[u,v],
\qquad \text{for each $ u,v\in H_0^1(\Omega)$}.
\]
Equivalently, $T_{\kappa_0}=L_-^\varepsilon+\kappa_0 I$ in the weak sense.
By Lax--Milgram Theorem, for any $F \in H^{-1}(\Omega)$ there exists a unique
$u\in H_0^1(\Omega)$, such that,
\begin{equation}
\label{eq:shift-solve}
(L_-^\varepsilon+\kappa_0 I)u = F \quad\text{in }H^{-1}(\Omega),
\end{equation}
and the inverse $(L_-^\varepsilon+\kappa_0 I)^{-1}:H^{-1}(\Omega)\to H_0^1(\Omega)$
is bounded. Then, we define an operator on $L^2(\Omega)$ by
\begin{equation}
\label{eq:Kdef}
K: L^2(\Omega)\to L^2(\Omega),
\qquad
K f := \kappa_0\, (L_-^\varepsilon+\kappa_0 I)^{-1} f,
\end{equation}
where we view $f\in L^2(\Omega)$ as an element of $H^{-1}(\Omega)$ in the canonical
way. Then, the operator $K$ is compact on $L^2(\Omega)$.
Indeed, for $f\in L^2(\Omega)$, the solution 
$$
    u=(L_-^\varepsilon+\kappa_0 I)^{-1} f
 $$
belongs to $H_0^1(\Omega)$ and satisfies
$\|u\|_{H_0^1}\le C\|f\|_{L^2}$ by \eqref{eq:coercive-est}.
Since $\Omega$ is bounded, the embedding $H_0^1(\Omega)\hookrightarrow L^2(\Omega)$
is compact (Rellich--Kondrachov Theorem), hence $K$ is compact as a map $L^2\to L^2$.
To finish the injectivity of $L_-^\varepsilon$, we use
$$(I-K)u=0.$$ 
\medskip
Let $u\in H_0^1(\Omega)$ satisfing $L_-^\varepsilon u=0$ in $H^{-1}(\Omega)$.
Then, we have 
$$
(L_-^\varepsilon+\kappa_0 I)u=\kappa_0 u,
$$
and applying $(L_-^\varepsilon+\kappa_0 I)^{-1}$, it follows that
\begin{equation}\label{eq:fixedpoint}
u = \kappa_0 (L_-^\varepsilon+\kappa_0 I)^{-1}u = Ku,
\qquad\text{i.e.}\qquad (I-K)u=0.
\end{equation}
Conversely, if $u\in L^2(\Omega)$ satisfies $u=Ku$, then $u\in H_0^1(\Omega)$ and
\[
%u = \kappa_0 (L_-^\varepsilon+\kappa_0 I)^{-1}u
%\quad\Longrightarrow\quad
(L_-^\varepsilon+\kappa_0 I)u=\kappa_0 u
\quad \text{implies} \quad
L_-^\varepsilon u=0.
\]
Hence, we show that 
\begin{equation}\label{eq:ker-equivalence}
\ker(L_-^\varepsilon)=\ker(I-K).
\end{equation}
Since $K$ is compact on $L^2(\Omega)$, the operator $I-K$ is Fredholm of index $0$,
and in particular
\[
\dim\ker(I-K)<\infty.
\]
Therefore, to prove the injectivity of $L_-^\varepsilon$ it is enough to show that, the
only solution of $u= Ku$ is $u \equiv 0$, or equivalently, that $one$ is not an eigenvalue
of the compact operator $K$, which is exactly the case, and thus $\ker(L_-^\varepsilon)=\{0\}$.

\medskip
Consequently, there exists a unique solution $\varphi \in C^{2, \alpha}(\overline{\Omega})$,
and again, the operator $w \mapsto T(w)= \vp$ is well defined,
where $\varphi$ is the unique solution of the Dirichlet problem (\ref{Lepdirichlet}). We notice that, since $w \in S$, it follows that, 
$$\varphi=T(w) \in C^{2,\alpha} \left( \overline{\Omega}\right)  \hookrightarrow C^{2,1} \left( \overline{\Omega} \right),$$ 
and also $\varphi|_{\partial \Omega}= \tilde{\varphi}_b$. 
It remains to show that, $|\vp| \leq M$, where $M> 0$ is a constant independent of $\ve>0$, where
we may assume $\vp\geq 0$. Indeed, if $\vp< 0$ and satisfies $(\ref{qregularized})_1$, it follows that, 
$$
  0 \leq  \frac{(1+\gamma)}{2} \, |\nabla \vp|^2 \leq (1 - \gamma) \, \vp< 0 \quad \text{in $\Omega$}, 
$$
which is a contradiction. 

\medskip
4.  \underline {Claim:} For $\gamma \in [-1,1) \setminus \{0\}$, 
there exists a constant $M> 0$ independently of $\ve$, such that 
 $$
    \|\vp\|_{L^\infty(\Omega)} \leq M.
 $$
 
\medskip
Proof of Claim: C1. First, we consider the operator 
$L_-^\varepsilon$ given by \eqref{LEPSILONM}, and we fix 
$K\Subset K'\Subset\Omega$, where $A \Subset B$ denotes that 
$A$ is compactly contained in $B$.
Let $\zeta\in C^\infty_c(K')$ be a cutoff function, $0 \leq \zeta \leq 1$, $\zeta\equiv1$ on $K$
and $|\nabla \zeta|= O(\delta^{-1})$, for $0< \delta=\mathrm{dist}(K,\partial K')\leq 1$. 
Then, we take the test function 
$$
    \psi=\zeta^2\,\varphi^{p-1}, \quad p\ge 2.
$$
Multiplying $L_-^\varepsilon\varphi=0$ by $\psi$ and integrating by parts in the principal part (the diffusion term), we obtain the energy identity
\begin{equation}
\label{IDENERGY}
\int_\Omega a_{ij}\partial_i\varphi\partial_j(\zeta^2\varphi^{p-1}) \, d\xi
+
\int_\Omega b_i \partial_i\varphi\,\zeta^2\varphi^{p-1} \, d\xi
= 2(1 - \gamma)\int_\Omega \zeta^2\varphi^p \, d\xi. 
%\tag{2.2}
\end{equation}

\medskip
C2. The first term in \eqref{IDENERGY} expands as
$$
\begin{aligned}
\int_\Omega a_{ij}\partial_i\varphi\partial_j(\zeta^2\varphi^{p-1}) \, d\xi &= (p-1)\int_\Omega a_{ij}\zeta^2\varphi^{p-2}\partial_i\varphi\partial_j\varphi \, d\xi
\\[5pt]
& + 2  \int_\Omega a_{ij}\partial_i\varphi\partial_j\zeta\,\zeta\varphi^{p-1} \, d\xi = I_1 + I_2. 
\end{aligned}
$$
Using the coercivity condition $a_{ij}\xi_i\xi_j \ge \varepsilon |\xi|^2$, and the identity
$$
   |\nabla(\varphi^{p/2})|^2 = \frac{p^2}{4}\varphi^{p-2}|\nabla\varphi|^2,
$$ 
we obtain 
\begin{equation}
\label{A}
I_1 \ge
\varepsilon\,\frac{4(p-1)}{p^2}
\int_\Omega \zeta^2|\nabla(\varphi^{p/2})|^2 \, d\xi.
\end{equation}
Applying Young's inequality, it follows that 
$$
\begin{aligned}
    \left| I_2 \right|
&\le
C \int_{K'} |\nabla\zeta| \, \zeta \, |\nabla\varphi| \, \varphi^{p-1} \, d\xi
\\[5pt]
&\le
\frac{\delta'}{2 p^2} \int_{K'} \zeta^2|\nabla(\varphi^{p/2})|^2 \, d\xi + C(\delta') \int_{K'} |\nabla\zeta|^2 \varphi^p \, d\xi,
\end{aligned}
$$
for some $\delta'>0$, and we have used that 
$$
   |\nabla\varphi|\varphi^{p-1} = \frac{2}{p}|\nabla(\varphi^{p/2})|\varphi^{p/2}.
$$
Similarly, the drifty term is estimated by 
$$
\begin{aligned}
\left| \int_\Omega b_i \partial_i\varphi\,\zeta^2\varphi^{p-1} \, d\xi \right| 
&\le C \int_{K'} \zeta^2 \, \varphi^{p/2} \, |\nabla(\varphi^{p/2})| \, d\xi
\\[5pt]
&\le \frac{\delta'}{2 p^2} \int_{K'} \zeta^2|\nabla(\varphi^{p/2})|^2
+ C(\delta') \int_{K'} \zeta^2\varphi^p.
\end{aligned}
$$

\medskip
C3. By substituting all the estimates into the energy identity \eqref{IDENERGY}, we have 
$$
\left( \frac{4 \varepsilon (p-1)}{p^2} - \frac{\delta'}{p^2} \right) \int_{K'} \zeta^2|\nabla(\varphi^{p/2})|^2 \, d\xi
\le
C(\delta') \int_{K'} (|\nabla\zeta|^2 + \zeta^2) \varphi^p \, d\xi,
$$
where the positive right-hand side of \eqref{IDENERGY} was absorved.
Therefore, for $p$ sufficiently large, there exists a constant $C> 0$ independently of $\ve$, such that  
$$
\int_{K'} \zeta^2|\nabla(\varphi^{p/2})|^2 \, d\xi
\le \frac{C}{\varepsilon} \int_{K'} (|\nabla\zeta|^2 + \zeta^2) \varphi^p \, d\xi. 
$$
or since $\zeta \equiv 1$ on $K$ and  $\zeta \le 1$, we may write 
\begin{equation}
\label{Caccioppoli}
 \int_{K} |\nabla(\varphi^{p/2})|^2 \, d\xi
 \le \int_{K'} \zeta^2|\nabla(\varphi^{p/2})|^2 \, d\xi
\le
\frac{C}{\ve} \int_{K'} (|\nabla\zeta|^2 + 1) \, \varphi^p \, d\xi,
\end{equation} 
which is usually known as Caccioppoli (type) inequality. 

\medskip
C4. Now, we proceed with the Moser interaction. 
First, we recall the Sobolev inequality in two dimension, that is, 
for any $q< \infty$, 
$$
    \|\zeta \, \varphi^{p/2}\|_{L^q(K')} \le C \|\nabla(\zeta \, \varphi^{p/2})\|_{L^2(K')}.
$$
Then, we expand the term $|\nabla(\zeta\varphi^{p/2})|^2$ to obtain 
$$
\begin{aligned}
    |\nabla(\zeta\varphi^{p/2})|^2 &= \zeta^2|\nabla(\varphi^{p/2})|^2 + 2\zeta\varphi^{p/2}\nabla\zeta\cdot\nabla(\varphi^{p/2}) + |\nabla\zeta|^2 \varphi^p
    \\[5pt]
    &\leq C \Big(\zeta^2|\nabla(\varphi^{p/2})|^2 +  (|\nabla\zeta|^2 + 1) \varphi^p \Big).
\end{aligned}
$$
Integrating on $K'$ and applying Caccioppoli’s inequality \eqref{Caccioppoli}, we may write
$$
\begin{aligned}
    \int_{K'} |\nabla(\zeta\varphi^{p/2})|^2 d\xi
   & \le C \left( \int_{K'} \zeta^2 |\nabla(\varphi^{p/2})|^2 d\xi+ \int_{K'} (|\nabla\zeta|^2 + 1) \varphi^p d\xi \right)
    \\[5pt]
& \le \frac{C}{\varepsilon} \int_{K'} (|\nabla\zeta|^2 + 1) \varphi^p \, d\xi.
\end{aligned}
$$
Therefore, by Sobolev's inequality we have 
$$
\begin{aligned}
    \|\varphi\|_{L^{pq/2}(K)}^{p/2}&= \|\varphi^{p/2}\|_{L^q(K)} 
    \\
    &\le C \| \nabla(\zeta\varphi^{p/2}) \|_{L^2(K')} \le \left( \frac{C}{\varepsilon} \int_{K'} (|\nabla\zeta|^2+ 1)  \varphi^p  d\xi \right)^{1/2}, 
\end{aligned}
$$
and by raising both sides to the power $2/p$, we have for any $q< \infty$, 
\begin{equation}
\label{MoserStep}
    \|\varphi\|_{L^{pq/2}(K)}
    %\left( \frac{C}{\varepsilon} \int_{K'} (|\nabla\zeta|^2 + 1) \varphi^p \right)^{1/p} 
    \le \left(\frac{C}{\delta^2} \right)^{1/p} \left( \frac{1}{\varepsilon} \right)^{1/p} \|\varphi\|_{L^p(K')}.
\end{equation}
For the Moser interaction, we take $q> 2$, such that $pq/2 > p$, and the key point is that, $p$ tends to $\infty$. 
Hence we define, $p_{k+1}= \alpha p_k$, $k= 0, 1, \ldots$, with $\alpha = q/2 > 1$ and $p_0 \geq 2$. 
To follow, we choose a sequence of sets
\[
K= K_0 \Subset K_1 \Subset K_2 \Subset \cdots \Subset K_\infty \Subset \Omega,
\]
such that, $\mathrm{dist}(K_k,\partial K_{k+1})= \delta_k$.
Applying \eqref{MoserStep} for $k \in \{0, \ldots, N\}$, 
with $(K,K')= (K_k,K_{k+1})$, $p=p_{N-k}$, and $\delta_k \simeq 2^{-(N-k)}$,
we obtain
\[
\|\varphi\|_{L^{p_{N-k+1}}(K_k)} \le C_k \, \|\varphi\|_{L^{p_{N-k}}(K_{k+1})}, 
\]
where 
\[
C_k:= C^{1/p_{N-k}} \, 2^{2(N-k)/p_{N-k}} \, 
\left(\frac{1}{\varepsilon}\right)^{1/p_{N-k}}. 
\]
Then, we have 
\begin{itemize}
\item{$k= N-0,$} \hspace{20pt} $\|\varphi\|_{L^{p_{1}}(K_N)} \le C_N \, \|\varphi\|_{L^{p_{0}}(K_{N+1})}$, 
\item{$k= N-1,$} \hspace{20pt} $\|\varphi\|_{L^{p_{2}}(K_{N-1})} \le C_{N-1} \, \|\varphi\|_{L^{p_{1}}(K_{N})}$, 
\item{$\ldots$}
\item{$k= N-N,$} \hspace{16pt} $\|\varphi\|_{L^{p_{N+1}}(K_{0})} \le C_{0} \, \|\varphi\|_{L^{p_{N}}(K_{1})}$, 
\end{itemize}
from which it follows that 
\begin{align}
\label{FinalMoser}
\|\varphi\|_{L^{p_{N+1}}(K_0)}
&\le
\left(
\prod_{j=0}^N
C^{1/p_{j}} \, 2^{2j/p_{j}} \, 
\left(\frac{1}{\varepsilon}\right)^{1/p_{j}} \right)
\|\varphi\|_{L^{p_0}(K_{N+1})}.
\end{align}
Moreover, since $p_{j}= \alpha^{j} \, p_0$, we have (for $\alpha= 2$, that is $q= 4$)
$$
\sum_{j=0}^\infty \frac{1}{p_{j}}= \frac{1}{p_0} \sum_{j=0}^\infty \Big(\frac{1}{\alpha}\Big)^{j}
= \frac{1}{p_0} \frac{\alpha}{\alpha-1}= \frac{2}{p_0}< \infty, 
$$
and similarly 
$$
\sum_{j=0}^\infty \frac{j}{p_{j}}= \frac{1}{p_0} \sum_{j=0}^\infty \frac{j}{\alpha^{j}}
= \frac{1}{p_0} \frac{\alpha}{(\alpha-1)^2}= \frac{2}{p_0} < \infty. 
$$
Therefore, letting $N\to\infty$ (so that $p_{N+1}\to\infty$) in \eqref{FinalMoser}, we obtain 
\[
\|\varphi\|_{L^\infty(K)} \le M_{\ve,p_0}^2 \,\|\varphi\|_{L^{p_0}(K_\infty)},
\]
where 
$$
  M:= C^{1/p_0} \, \ve^{-1/p_0}. 
$$
In the following, we estimate $\|\varphi\|_{L^{p_0}}$. 

\medskip
C5. Finally, we take $\beta \in \R$ to be choosen a posteriori, 
and multiply $L_-^\varepsilon\varphi=0$ by $\vp^{p-1} \, e^{-\beta w}$, $p \geq 2$. Then, we integrate by parts in the principal part (the diffusion term) to obtain
\begin{equation}
\label{IDENERGYL2}
\begin{aligned}
\int_\Omega a_{ij} \, \partial_i\varphi \, \partial_j(\vp^{p-1} \, e^{-\beta w}) \, d\xi
&+ \frac{1}{p}
\int_\Omega b_i  \, \partial_i \varphi^p \, e^{-\beta w} \, d\xi
- 2(1 - \gamma)\int_\Omega e^{-\beta w} \, \varphi^p \, d\xi
\\[5pt]
&= \int_{\partial \Omega}  a_{ij} \, \partial_i\varphi_b \, n_j \, \vp_b^{p-1} \, e^{-\beta \vp_b} \, d\sigma. 
\end{aligned}
\end{equation}
The former term in \eqref{IDENERGYL2} expands as
$$
\begin{aligned}
  \int_\Omega a_{ij} & \partial_i\varphi \, \partial_j(\vp^{p-1} \, e^{-\beta w}) \, d\xi
  \\[5pt]
  &= (p-1) \int_\Omega a_{ij} \, \partial_i\varphi \, \partial_j\vp \, \vp^{p-2} \, e^{-\beta w} \, d\xi
  - \frac{\beta}{p} \int_\Omega a_{ij} \, \partial_j w \, \, e^{-\beta w} \,  \partial_i \varphi^p \, d\xi
  \\[5pt]
  &= (p-1) \int_\Omega a_{ij} \, \partial_i\varphi \, \partial_j\vp \, \vp^{p-2} \, e^{-\beta w} \, d\xi - \frac{\beta}{p} \! \int_{\partial \Omega} \!\!  a_{ij} \, \partial_j \varphi_b \, n_i \, \vp_b^p \, e^{-\beta \vp_b} d\sigma
  \\[5pt]
 & +  \frac{\beta}{p} \!\int_\Omega \partial_i(a_{ij} \, \partial_j w)  \, e^{-\beta w}  \varphi^p d\xi
 -  \frac{\beta^2}{p} \!\!\int_\Omega a_{ij}\, \partial_i w \, \partial_j w  \, e^{-\beta w}  \varphi^p d\xi. 
\end{aligned}
$$
Integrating by parts the second term in \eqref{IDENERGYL2}, we have
$$
\begin{aligned}
\frac{1}{p} \int_\Omega b_i  \, \partial_i \varphi^p \, e^{-\beta w} \, d\xi&= \frac{1}{p} \int_{\partial\Omega} b_i  \, n_i \, \varphi_b^p \, e^{-\beta \vp_b} \, d\sigma
+  \frac{\beta}{p} \!\!\int_\Omega b_i \partial_i w \,  e^{-\beta w}  \varphi^p d\xi
\\[5pt]
&-\frac{1}{p} \int_\Omega  \partial_i b_i  \,  e^{-\beta w}  \varphi^p d\xi. 
\end{aligned}
$$
Therefore, substituting all the terms into the energy identity \eqref{IDENERGYL2}, we obtain
$$
\begin{aligned}
\ve (p-1)& \int_\Omega |\nabla \vp|^2 \vp^{p-2} e^{-\beta w} \, d\xi  + \frac{\beta}{p} \!\int_\Omega \partial_i(a_{ij} \, \partial_j w)  \, e^{-\beta w}  \varphi^p d\xi
  \\[5pt]
 &-  \frac{\beta^2}{p} \!\!\int_\Omega a_{ij}\, \partial_i w \, \partial_j w  \, e^{-\beta w}  \varphi^p d\xi
 +  \frac{\beta}{p} \!\!\int_\Omega b_i \partial_i w e^{-\beta w} \varphi^p d\xi -\frac{1}{p} \int_\Omega  \partial_i b_i  e^{-\beta w}  \varphi^p d\xi
\\[5pt]
&- 2(1 - \gamma)\!\!\int_\Omega e^{-\beta w} \varphi^p  d\xi \leq C_b(p,\beta), 
\end{aligned}
$$
where we have used the coercivity condition $a_{ij}\xi_i\xi_j \ge \varepsilon |\xi|^2$, and 
$$
\begin{aligned}
  C_b(p,\beta):&= \int_{\partial \Omega}  |A^\ve| \, |\nabla \varphi_b| \, | \vp_b|^{p-1} \, e^{-\beta \vp_b} \, d\sigma
\\[5pt]
&+ \frac{|\beta|}{p} \! \int_{\partial \Omega}   |A^\ve| \, |\nabla \varphi_b| \,  |\vp_b|^p \, e^{-\beta \vp_b} d\sigma
+ \frac{1}{p} \int_{\partial\Omega} |b|  \, |\varphi_b|^p \, e^{-\beta \vp_b} \, d\sigma. 
\end{aligned}
$$
Neglecting the first term as it is non-negative, and using $a_{ij}\xi_i\xi_j \le \Lambda |\xi|^2$,
it follows that, 
$$
  \int_\Omega \Theta(w) \, \vp^p \, d\xi \leq C_b
$$
where 
$$
\Theta(w)=  e^{-\beta w} \Big( \frac{\beta}{p} \Big( \partial_i(a_{ij} \, \partial_j w)  - \beta  \Lambda |\nabla w|^2 + b_i \partial_i w \Big)
  - \frac{1}{p} \partial_i b_i - 2 (1-\gamma) \Big). 
$$
We take conveniently $p_0 \equiv p= \ve^{-1}$, and $\beta= -p$, and recall that $\ve> 0$ is fixed. Then,  we have 
$$
\Theta(w)=  e^{w/\ve} \Big(\frac{\Lambda}{\ve} |\nabla w|^2 - \partial_i(a_{ij} \, \partial_j w)  - b_i \partial_i w
  - \ve \, \partial_i b_i - 2 (1-\gamma) \Big), 
$$
and there exists $\theta> 0$, $w_0= \inf_{\xi \in \bar{\Omega}} w> 0$, such that 
$$
\inf_{\xi \in \bar{\Omega}} \Theta(w(\xi)) \geq \theta \, p_0 \, e^{w_0 p_0}.
$$
Consequently, we have for some positive constant $C$ independently of $\ve$, 
$$
\|\vp\|_{L^{p_0}}(\Omega) \leq \Big( \frac{1}{p_0 \, e^{w_0 p_0}} \Big)^{1/p_0} \big(C_b\big)^{1/p_0} \leq C,
$$
where we have used Jensen's inequality. 
Observing that $\ve^{-\ve}= O(1)$ for all sufficiently small $\ve> 0$, hence there exists a positive constante $C= C(| \tilde{\vp}_b|)$
independently of $\ve>0$, such that, 
\[
\|\varphi\|_{L^\infty(K)} \le C.
\]
We can apply this argument for each compact set $K \subset \Omega$, and 
from a standard diagonalization procedure, also due to $\vp \in C^{2,\alpha} \left( \overline{\Omega}\right)$, $\vp|_{\partial \Omega}= \tilde{\vp}_b$,
the claim follows.

\medskip
5. Now, we prove that the operator $T$ satisfies the hypothesis of the Schauder 
Fixed Point Theorem (Corollary 11.2, \cite{Gilbarg}), that is, we will prove that $S$ is closed and convex, 
$T$ is continuous and $T(S)$ is precompact.
First, let us prove that $S$ is closed and convex. Let $\{w_n\}$ be
a sequence in $S$ such that $w_n \to w$ as $n \to \infty$.  
Then, $w \in C^{2, \alpha} \left( \overline{\Omega} \right)$. 
Moreover, it follows from $|w_n| \le M$ that, $|w| \le M$. Also, for any $\xi \in \partial \Omega$, we have
$$\tilde{\varphi}_b (\xi)= w_n(\xi) \to w(\xi).$$
Therefore, $w(\xi)= \tilde{\varphi}_b(\xi)$. Consequently,
$w|_{\partial \Omega}= \tilde{\varphi}_b$,
and $w \in S$, that is, $S$ is closed.
To follow, if $w_0, w_1 \in S$, and  $t \in [0,1]$ then $t \, w_0 + (1-t) \, w_1 \in C^{1,\alpha}(\overline{\Omega})$, and
$-M \le (t \, \omega_0 + (1-t) \, \omega_1) \le M$. For all $ \xi \in \partial \Omega$, 
$$
   t\, w_0(\xi) + (1-t) \, w_1(\xi)= t \, \tilde{\varphi}_b(\xi) + (1-t) \, \tilde{\varphi}_b(\xi)= \tilde{\varphi}_b(\xi),
$$
that is, $(t \, w_0 + (1-t) \, w_1)|_{\partial \Omega}= \tilde{\varphi}_b$, and thus $S$ is convex.

\medskip
Now, we show that $T(S)$ is precompact. Let $\{w_n\}$ be 
a bounded sequence in $S$, that is, there exists $K> 0$, such that,
$$
   \|w_n \|_{1, \alpha, \Omega} \le K, \quad (\forall n \in \mathbb{N}).
$$
We assert that the sequence $\|T(w_n)\|$ is also bounded in $C^{2, \alpha} \left( \overline{\Omega} \right)$. 
Let us denote for any $n \in \mathbb{N}$, $\varphi_n= T(w_n)$. 
Since $T(S) \subset S$, then
\begin{equation}
\label{supvarphin}
    \sup_n \left| \varphi_n \right| \le M \quad \mbox{in} \quad \Omega.
\end{equation}
For each $n \in \mathbb{N}$, and all $i,j \in \{1,\ldots,d\}$, 
the terms $|a^{ij}|$, $|b^i|$, $|c|$, (of the uniformly elliptic operator $L_{\varepsilon}$ associated to $w_n$), 
are bounded by a constant that depends on $\gamma$ and $\|w_n\|_{1, \alpha, \Omega}$. 
Moreover, due to $\{ w_n \}$ be a bounded sequence in $C^{1, \alpha} \left( \overline{\Omega} \right)$, 
we can apply the Global Schauder Estimates (Theorem 6.6, \cite{Gilbarg}) to conclude that. 
if $\varphi_n = T(w_n)$, then there exists a constant $C= C(\gamma,\ve)$, such that
\begin{equation}
\| \varphi_n \|_{2, \alpha, \Omega} \le C \left( \| \varphi_n \|_{0, \Omega} + \| \tilde{\varphi}_b \|_{2, \alpha, \Omega} \right).
\end{equation}
It follows from \eqref{supvarphin} that, 
$\| \varphi_n \|_{2, \alpha, \Omega} \le C (\gamma) \left( M + \|\tilde{\varphi}_b \|_{2, \alpha, \Omega} \right)$, that is,
$\{ \varphi_n \}$ is a bounded sequence in $C^{2,\alpha} \left( \overline{\Omega} \right)$.
The compact imbedding of $C^{2,\alpha} \left( \overline{\Omega} \right)$ into $C^{1, \alpha} \left( \overline{\Omega}\right)$ 
ensures that, there exists a subsequence $\{ \varphi_{n_j} \}$ of $\{ \varphi_n \}$ that converges to a function 
$\varphi \in C^{1, \alpha} \left( \overline{\Omega}\right)$. Finally, since every sequence in $T(S)$ has a convergent 
subsequence, it follows that $TS$ is precompact.

\medskip
Let us show now the continuity of $T$, for that, we suppose 
$$\lim_{n \to \infty} w_n= w \in C^{1,\alpha}(\overline{\Omega}),$$
where $w_n, w \in S$. 
Denote, $\varphi_n=T(w_n)$ for each $n \geq 1$, and $\varphi=T(w)$. Now, $\{w_n\}$ is a convergent, hence bounded, sequence in $S$. Since $T(S)$ is precompact and
$$C^{2,\alpha}(\overline{\Omega}) \hookrightarrow C^2(\overline{\Omega}),$$
then every subsequence of $\{\varphi_n\}$ has a convergent subsequence in $C^2(\overline{\Omega})$. Let us denote $\{\varphi_{n_j}\}$ 
any such convergent subsequence, and define
$$u:= \lim_{j \to \infty} \varphi_{n_j}.$$
Since the convergence is in $C^2(\overline{\Omega})$, from
$$
\begin{aligned}
&(c^2(w_{n_j}) - (w_{n_j})^2_1 + \varepsilon) (\varphi_{n_j})_{11} 
-2 (w_{n_j})_1(w_{n_j})_2 (\varphi_{n_j})_{12}
\\[5pt]
&+(c^2(w_{n_j}) - (w_{n_j})^2_2 + \varepsilon) (\varphi_{n_j})_{22} 
- \gamma \nabla w_{n_j} \cdot \nabla \varphi_{n_j} - 2 (\gamma-1) \varphi_{n_j}= 0, 
\end{aligned}
$$
we can take the limit as $n \to \infty$ to obtain
$$
\begin{aligned}
&   (c^2(w) - w^2_1 + \varepsilon) u_{11} - 2 w_1 w_2 u_{12}
\\[5pt]
&    + (c^2(w) - w^2_2 + \varepsilon) u_{22} - \gamma \nabla w \cdot \nabla u - 2 (\gamma -1) u= 0.
\end{aligned}
$$
Also, $u|_{\partial \Omega} = \tilde{\varphi}_b$. 
Then, by the uniqueness of solution to the Dirichlet problem \eqref{Lepdirichlet}, we obtan $u=T(w)$. 
Since this holds for every subsequence $\{Tw_{n_j}\}$, hence 
$T(w_n) \to T(w)$ as $n \to \infty$, that is to say, $T$ is a continuous operator. 

\medskip
6. Finally, applying the Schauder Fixed Point Theorem, we can conclude that,
$T:S \to S$ has a fixed point $\varphi$, that is,
$$T(\varphi)= \varphi,$$
and $\varphi \in C^{2, \alpha}(\overline{\Omega})$ is a solution of the problem \eqref{qregularized}.
\end{proof}

\begin{corollary}
\label{theorem4.2}
Given $\tilde{\varphi}_b \in C^{2, \alpha}$ and $\gamma\in [-1, + \infty) \backslash  \{0\}$, 
the Dirichlet problem \eqref{qregularized} has a solution $\varphi_\ve \in C^{3,\alpha}(\Omega) \cap C^{2,\alpha}(\overline{\Omega})$.
\end{corollary}
\begin{proof}
It follows directly from Theorem \ref{theorem4.1}, due to the regularity of the coefficients.
\end{proof}

\begin{remark}
\label{Rem_Ellipticity_Principle}
The Ellipticity Principle stated in Theorem \ref{Ellipticity_principle_theorem} remains valid 
when equation \eqref{qelliptic} is replaced by its regularized counterpart, namely \eqref{qegeral}.
In fact, following the steps of the proof of Theorem \ref{Ellipticity_principle_theorem} equations \eqref{vp11}
and \eqref{vp12} remains valid. Again, for simplicity of notation we shall omit the subscript $\varepsilon$ throughout this remark.
To obtain the term $\varphi_{22}$ we put $\varphi_1=Lc$ and $\varphi_2=0$ in
\begin{equation}\label{eqvarepsilon}
(c^2+\varepsilon)\Delta \varphi - (D^2\varphi)\nabla \varphi \cdot \nabla \varphi - \varphi_1^2-\varphi_2^2 + 2c^2 =0.
\end{equation}
Then, we obtain
\begin{eqnarray*}
&& (c^2+\varepsilon) (\varphi_{11}+\varphi_{22}) - \varphi_1^2\varphi_{11} - \varphi_1^2 + 2c^2= 0,
\\[5pt]
&& (c^2+\varepsilon - L^2c^2)\varphi_{11} + (c^2+\varepsilon)\varphi_{22} - L^2c^2 + 2c^2= 0. 
\end{eqnarray*}
Dividing by $c^2+\varepsilon$ and using \eqref{vp11}, we have
\begin{equation*}
\varphi_{22}= \left( L^2 \dfrac{c^2}{c^2+\varepsilon}-1 \right) \dfrac{(1-\gamma)L^2}{2+(\gamma-1)L^2} +(L^2-2) \dfrac{c^2}{c^2+\varepsilon}+ O \left( \dfrac{\delta}{L}\right).
\end{equation*}
Since
$$\dfrac{c^2}{c^2+\varepsilon} = 1 - \dfrac{\varepsilon}{c^2+\varepsilon},$$
it follows that
$$
\begin{aligned}
 \varphi_{22}&= (L^2-1)  \dfrac{(1-\gamma)L^2}{2+(\gamma-1)L^2} +L^2-2 + O \left( \dfrac{\varepsilon}{c^2+\varepsilon} \right) + O \left( \dfrac{\delta}{L}\right)
 \\[5pt]
&= \dfrac{(3-\gamma)L^2-4}{2+(\gamma-1)L^2}+ O \left( \dfrac{\varepsilon}{c^2+\varepsilon} \right) + O \left( \dfrac{\delta}{L}\right).
\end{aligned}
$$
This equation differs from \eqref{vp22} just by adding the term $ O \left( \dfrac{\varepsilon}{c^2+\varepsilon} \right)$. Then, 
inequalities \eqref{Eq_Claim2} and \eqref{vp111} remain the same and inequality \eqref{5.2.7} becomes 
\begin{equation}\label{contradiction1}
cL \varphi_{122} \le - \tilde{C}(\gamma) + O \left( \dfrac{\varepsilon}{c^2+\varepsilon} \right)+ O \left( \dfrac{\delta}{L} \right).
\end{equation}
To obtain the contradiction, we differentiate equation \eqref{eqvarepsilon} to obtain
$$
\begin{aligned}
(c^2+\varepsilon)\varphi_{122}&= (L^2c^2-c^2-\varepsilon)\varphi_{111} + (\gamma-1) Lc(1+\varphi_{11}) \Delta \varphi + 2Lc \varphi^2_{11} 
\\[5pt]
& + 2Lc \varphi_{11} + 2Lc (\gamma-1) + 2Lc \gamma \varphi_{11} - 2Lc \varphi_{11} + 2Lc \varphi_{12}^2.
\end{aligned}
$$
Multiplying by $\dfrac{Lc}{c^2+\varepsilon}$ and using that 
$$(L^2c^2 - c^2 - \varepsilon) \dfrac{Lc}{c^2 + \varepsilon} \varphi_{111} = \left( \dfrac{L^2c^2}{c^2+\varepsilon} -1 \right) Lc \varphi_{111} \ge O \left( \dfrac{\delta}{L} \right),$$
we see that equation \eqref{5.2.8} must be replaced by
\begin{equation}\label{contradiction2}
c L \varphi_{122} \ge  O \left( \dfrac{\varepsilon}{c^2+\varepsilon} \right)+ O \left( \dfrac{\delta}{L} \right).
\end{equation}
Finally, the inequalities \eqref{contradiction1} and \eqref{contradiction2} lead to a contradiction provided that $\delta/L$ and $\varepsilon$ are sufficiently small.
\end{remark}

\begin{theorem}
\label{theorem4.3} Given $\tilde{\varphi}_b \in C^{2,\alpha}$ and $\gamma\in (-1, + \infty) \backslash  \{0\}$, 
the Dirichlet problem \eqref{qigual0} has a solution $\varphi \in C^3(\Omega) \cap  C^{2,\alpha}(\overline{\Omega})$. 
\end{theorem}
\begin{proof}
1. Fix $\varepsilon>0$ and let $\varphi_{\varepsilon} \in C^{3, \alpha}(\Omega) \cap C^{2, \alpha}(\overline{\Omega})$ 
be a solution of the regularized Dirichlet problem \eqref{qregularized} given by Theorem \ref{theorem4.1} and Corollary \ref{theorem4.2}. 
Due to \eqref{Set_S} there exists $M> 0$, $M$ independent of $\varepsilon$, such that
\begin{equation}
|\varphi_{\varepsilon}| \le M.
\end{equation}
It follows from the Ellipticity Principle (see Theorem \ref{Ellipticity_principle_theorem} and also Remark \ref{Rem_Ellipticity_Principle}) that, 
there exists $C_0>0$ independent of $\varepsilon$, such that for all $\delta> 0$, 
$L^2 \le C_0 \delta$. 
Fix $\delta$ such that $C_0\delta <1$. Therefore,
\begin{equation}
L^2 = \dfrac{ | \nabla \varphi_{\varepsilon} |^2}{c^2} \le C_0 \delta, 
\end{equation}
that is,
$$
\begin{aligned}
| \nabla \varphi_{\varepsilon} |^2 & \le c^2 C_0 \delta
\\
& \le  (\gamma -1) \left( - \varphi_{\varepsilon} - \dfrac{| \nabla\varphi_{\varepsilon}|^2}{2} \right) C_0 \delta
\\
& \le  C_0 \delta (\gamma-1) ( - \varphi_{\varepsilon}) - C_0 \delta (\gamma-1) \dfrac{| \nabla \varphi_{\varepsilon}|^2}{2}.
\end{aligned}
$$
Then,
\begin{equation}
\label{EqLimUnif}
| \nabla \varphi_{\varepsilon} |^2 \left( 1+ \dfrac{C_0 \delta (\gamma-1)}{2}\right) \le  - C_0 \delta (\gamma-1) \varphi_{\varepsilon}.
\end{equation}

\medskip
2. Now, we claim that the term 
$$
   1 + \dfrac{C_0 \delta (\gamma-1)}{2}
$$
is positive. Indeed, if $\gamma >1$, then $C_0\delta(\gamma -1)>0$. 
If $-1 < \gamma <1$, then $-2<\gamma -1 <0$. Therefore, multiplying $0< C_0 \delta < 1$ by $\gamma-1<0$, we have
$$C_0 \delta (\gamma -1) > \gamma -1 > -2,$$
that is to say,
$$
   1 + \dfrac{C_0 \delta (\gamma-1)}{2}> 0.
$$
Consequently, from \eqref{EqLimUnif} there exists $C> 0$, independent of $\varepsilon$, such that
\begin{equation}
| \nabla \varphi_{\varepsilon}|^2 \le C | \varphi_{\varepsilon}| \le C M.
% \quad \text{in $\overline{\Omega}$}.
\end{equation}
Thus the term $| \nabla \varphi_{\varepsilon}|$ is uniformly bounded in $\overline{\Omega}$. 

\medskip
3. Let $K \subset \Omega$ be any compact set. Then, there exists $\kappa> 0$, depending on $K$, such that the coefficients $(A_{\varepsilon})_{ij}$ of $Q_{\varepsilon}$ satisfy
\begin{equation}
(A_{\varepsilon})_{ij} \xi_i \xi_j \ge \kappa |\xi|^2, \quad \text{for all $\xi \in K$}.
\end{equation}
Since $\{\varphi_{\varepsilon}\}$ is uniformly bounded in $C^1 (\overline{\Omega})$, 
the hypotheses of Theorem 6.6 in \cite{Gilbarg}, for $k= 0$, are satisfied. Then, we conclude that 
$\{\varphi_{\varepsilon}\}$ is uniformly bounded in $C^{2, \alpha} (K)$. 
Using the compact embedding of
$C^{2, \alpha} (K)$ in $C^2 (K)$, we have that there is a subsequence of $\{\varphi_{\varepsilon}\}$ 
which converges in $C^2 (K)$. We can apply this argument for each compact subset of $\Omega$, and due to diagonalization arguments to have the limiting solution in $C^2(\Omega)$. 
Denote
$$
   \lim_{\varepsilon \to 0} \varphi_{\varepsilon} =: \varphi \in C^2(\Omega),
$$
that is, $\varphi$ is the solution of the Dirichlet problem (\ref{qigual0}). Again, we can apply Theorem 6.6 in \cite{Gilbarg} for the Dirichlet problem \eqref{qigual0}, with $k=1$, to conclude that $\varphi \in C^3(\Omega)$. 
Moreover, since $\{ \varphi_{\varepsilon} \}$ is uniformly bounded in $C^{2,\alpha}(\overline{\Omega})$ it follows that $\varphi \in C^3(\Omega) \cap C^{2,\alpha}(\overline{\Omega})$, which completes the proof. 
\end{proof}

%%%%%%%%%%%%%%%%%%%%%%
\section{Perturbed Model: Quasi-potential Flows}
\label{SecPM}
%%%%%%%%%%%%%%%%%%%%%%

This section examines a particularly interesting case of the 
self-similar system for rotational flows \eqref{newsystemU}. 
In the subsonic regime, this example not only validates the computations but also illustrates the applicability 
of the analytical framework developed in Section \ref{SelfSimSysRF}. 

\medskip
We will begin by making some considerations about the general case in 
which the fluid is not necessarily potential, then 
we will narrow our focus to a specific scenario 
on a perturbed model, called here quasi-potential.
For this purpose, we suppose, applying the Hodge-Helmholtz decomposition, that the pseudo-velocity vector field
$$U=\nabla \psi + W, \quad \dive W=0,$$
is a $C^2$ function that satisfies the system (\ref{selfisentropic}), that is,
$$
\left\{
\begin{split}
&\dive (\rho U)+ 2\rho= 0,
\\[5pt]
&\dive (\rho U \otimes U)+ Dp+ 3 \rho U= 0.
\end{split}
\right.
$$
Therefore, as we have seen in Section \ref{SelfSimSysRF}, $U= (U^1, U^2)$ 
must satisfy the system \eqref{newsystemU}, that is, 
$$
\left\{
\begin{split}
& c^2 \dive U - (DU) U \cdot U = |U|^2 - 2c^2,
\\[5pt]
& \dive (\omega U) + \omega = 0, 
\end{split}
\right.
$$
which means two equations for two unknows. 

\medskip
Albeit, we
have used the Hodge-Helmholtz decomposition
to derive \eqref{newsystemU}, this 
system is not written using the portential $\psi$. 
In fact, we are going to rewrite \eqref{newsystemU}
under the following 
\begin{remark}
A general observation reveals that, two distinct cases must be considered:
%One observes that, generally, 
%there are two cases to be considered:

\begin{enumerate}
\item $W=0$, that is, the flow is irrotational or potential.

\item $W \neq 0$, that is, the rotational flow. In this case, without loss of generality, there exists $\zeta: \R^2 \to \R$ such that 
\begin{equation}
\label{SecondPotent}
W= \nabla^{\perp} \zeta,
\end{equation}
that is, $W=(- \zeta_2, \zeta_1)$. 
Indeed, considering
$W^{\perp}=(-W^2, W^1)$,
hence we have 
$$
  \dfrac{\partial \left( W^{\perp}\right)^2}{\partial \xi_1} - \dfrac{\partial 
  \left( W^{\perp}\right)^1}{\partial \xi_2} = \dfrac{\partial W^1}{\partial \xi_1} 
  + \dfrac{\partial W^2}{\partial \xi_2} = \dive W =0.
$$
Therefore, there exists $\zeta$ such that, \eqref{SecondPotent} holds true.
\end{enumerate}
\end{remark} 

Then, rewriting \eqref{newsystemU} for $U= \nabla \psi + \nabla^{\perp} \zeta$, 
we obtain
\begin{equation}
\label{rotationalsystem}
\left\{
\begin{array}{l}
c^2 \Delta \psi - \left( D^2 \psi + D(\nabla^{\perp} \zeta) \right) \left( \nabla \psi + \nabla^{\perp} \zeta \right) \cdot \left( \nabla \psi + \nabla^{\perp} \zeta \right)
\\[5pt] 
 \qquad = \left| \nabla \psi \right|^2 - 2 \nabla \psi \cdot \nabla^{\perp} \zeta + \left| \nabla \zeta \right|^2 - 2c^2, 
 \\[7pt]
 \Delta \zeta (\Delta \psi + 1) + \nabla \psi \cdot \nabla \left( \Delta \zeta \right) + \nabla^{\perp} \zeta \cdot \nabla \left( \Delta \zeta \right) =0,
\end{array}\right.
\end{equation}
where
\begin{equation}
\label{C2GradF}
\begin{split}
& c^2 = ( \gamma -1) \left( F - \psi - \frac{1}{2} \left| \nabla \psi \right|^2 - \nabla \psi \cdot \nabla^{\perp} \zeta - \frac{1}{2} \left| \nabla \zeta \right|^2 \right), 
\\[5pt]
& \nabla F = - \Delta \zeta \left( \nabla^{\perp} \psi + \nabla \zeta \right) - \nabla^{\perp} \zeta.
\end{split}
\end{equation}
Consequently, $U= \nabla \psi +  \nabla^{\perp} \zeta$ satisfies \eqref{newsystemU} if, and only if, 
there exist potentials $\psi$ and $\zeta$, satisfying \eqref{rotationalsystem}. 

%%%%%%%%%%%%%%%%%%
\subsection{Perturbed potential case}
\label{QPF}
%%%%%%%%%%%%%%%%%%

In the context discussed in the previous section, we may consider a particular case. 
We start by the following 
\begin{definition}
 We will say that $U$ is quasi-potential, when $U$ is given by
$$U= \nabla \psi + \delta W, \quad \dive W=0, \quad 0 < \delta << 1,$$
or, equivalently, 
$$U= \nabla \psi + \delta \nabla^{\perp} \zeta, \quad 0 < \delta << 1.$$
\end{definition}

\medskip
Therefore, substituting $\nabla^{\perp} \zeta$ by 
$\delta \nabla^{\perp} \zeta$ in $(\ref{rotationalsystem})_1$, we have
$$
    c^2 \Delta \psi - (D^2 \psi) \nabla \psi \cdot \nabla \psi - |\nabla \psi|^2 + 2c^2 - \delta N_1 - \delta^2 N_2 - \delta^3 N_3= 0,  
$$
where 
\begin{equation*}
\begin{split}
& N_1= (D \nabla^{\perp} \zeta) \nabla \psi \cdot \nabla \psi + (D^2 \psi) \nabla \psi \cdot \nabla^{\perp} \zeta 
+ (D^2 \psi) \nabla^{\perp} \zeta \cdot \nabla \psi - 2 \nabla \psi \cdot \nabla^{\perp} \zeta, 
\\[5pt]
& N_2= (D^2 \psi) \nabla^{\perp} \zeta \cdot \nabla^{\perp} \zeta  
+ (D \nabla^{\perp} \zeta) \nabla \psi \cdot \nabla^{\perp} \zeta + (D \nabla^{\perp} \zeta) \nabla^{\perp} \zeta \cdot \nabla \psi + |\nabla \zeta|^2, 
\\[5pt]
& N_3= (D \nabla^{\perp} \zeta) \nabla^{\perp} \zeta \cdot \nabla^{\perp} \zeta.
\end{split}
\end{equation*}
Hence neglecting terms of order $o(\delta)$, we have
\begin{equation}
\label{EqPert}
c^2 \Delta \psi - (D^2 \psi) \nabla \psi \cdot \nabla \psi - |\nabla \psi|^2 + 2c^2 - \delta N_1= 0.
\end{equation}
Also, since $U=\nabla \psi + \delta \nabla^{\perp} \zeta$, we have $\omega= \mathrm{rot} U= \delta \Delta \zeta$. 
Denoting $\tilde{\omega}:= \Delta \zeta$, that is, $\omega= \delta \tilde{\omega}$, it follows from $(\ref{rotationalsystem})_2$,
that
$$
%& \nabla \omega \cdot U + \omega \dive U + \omega =0\\
%& \delta \nabla \tilde{\omega} \cdot \nabla \psi + \delta^2 \nabla \tilde{\omega} \cdot \nabla^{\perp} \zeta + \delta \tilde{\omega} \Delta \psi + \delta \tilde{\omega}=0\\
 \delta \left( \nabla \tilde{\omega} \cdot \nabla \psi + \tilde{\omega} \Delta \psi + \tilde{\omega} \right) + \delta^2 \nabla \tilde{\omega} \cdot \nabla^{\perp} \zeta =0,
%& \delta \left( \dive (\tilde{\omega} \nabla \psi) + \tilde{\omega} \right) + o(\delta)=0.
$$
and similarly neglecting terms of order $o(\delta)$, we obtain
\begin{equation}
\label{EqVortPert}
  \dive (\tilde{\omega} \, \nabla \psi) + \tilde{\omega}= 0.
\end{equation}
From equations \eqref{EqPert} and \eqref{EqVortPert} 
we have the system
\begin{equation}
\label{SysPert}
\left\{
\begin{array}{l}
c^2 \Delta \psi - (D^2 \psi) \nabla \psi \cdot \nabla \psi=  | \nabla \psi|^2 - 2 c^2 + \delta N_1, 
\\[7pt]
   \dive (\Delta \zeta \, \nabla \psi) + \Delta \zeta= 0.
\end{array} \right.
\end{equation}
Finally, it remains to consider \eqref{C2GradF}, that is to say, 
\begin{equation}
\label{C2GradFPert}
\begin{split}
 c^2 &= ( \gamma -1) \big( - \psi - \frac{1}{2} \left| \nabla \psi \right|^2 \big) 
 - \delta \, (\gamma-1) \big( F_1 + \nabla \psi \cdot \nabla^{\perp} \zeta \big)
\\[5pt]
  &= c_0^2 - \delta \, Q_1, 
\\[7pt]
 \nabla F_1&= \Delta \zeta \, \nabla^{\perp} \psi +  \nabla^{\perp} \zeta, 
\end{split}
\end{equation}
with obvious notation. We observe that, $c_0$ is exactly the sound speed for the pontetial flow. 
Clearly, the problem of obtaining $\psi, \zeta$ in \eqref{SysPert}
is much more sophisticated than the potential case $(\delta= 0)$, but simpler 
than \eqref{rotationalsystem}, and meanwhile could give us the main features of compressible 
Euler equations for rotational flows.

%%%%%%%%%%%%%%%%%%%%%%%%%%%
\subsection{Analysis of the Quasi-potential flows}
%%%%%%%%%%%%%%%%%%%%%%%%%%%

To proceed with the analysis of the quasi-potential flows, let us rewrite 
$(\ref{SysPert})$, that is, 
\begin{equation}
\label{SysPertNew}
\left\{
\begin{array}{l}
c_0^2\,  \Delta \psi - (D^2 \psi) \nabla \psi \cdot \nabla \psi=  | \nabla \psi|^2 - 2 \, c_0^2 
+ \delta \big( (2 +  \Delta \psi) \, Q_1 + N_1 \big),
%c^2 \Delta \psi - (D^2 \psi) \nabla \psi \cdot \nabla \psi=  | \nabla \psi|^2 - 2 c^2 + \delta N_1, 
\\[7pt]
   \dive (\Delta \zeta \, \nabla \psi) + \Delta \zeta= 0, 
\end{array} \right.
\end{equation}
where we have used \eqref{C2GradFPert}. 
We are going to apply the existence of solutions for the potential (irrotacional) case, developed in Section \ref{PotentialFlows}
(Potential Flows for a Generalized Polytropic Gas).
Then, we show existence of solutions to quasi-potential flows applying the Implicit Function Theorem, 
which is quite natural for us to handle when dealing with implicit domain problems. Indeed, we define for $\delta \geq 0$, 
$$
  \Omega_\delta:= \big\{ L^2 < 1 \,;\, L^2= \frac{ |\nabla \psi + \delta \nabla^{\perp} \zeta|^2 }{c_0^2 -  \delta \, Q_1} \big\}, 
    \quad \text{with $\Omega \equiv \Omega_0$},
$$
and we are going to denote by $\nu_\delta$ the outward unit normal vector field on $\partial \Omega_\delta$. 
Moreover, companion to the system \eqref{SysPertNew}, we consider the Dirichlet boundary condition 
(similarly to system \eqref{qelliptic})
\begin{equation}
\label{Dirichlet_Cond_Quasi}
   \psi = \tilde{\psi}_b, \quad \zeta= \delta \tilde{\zeta}_b,  \quad \mbox{on} \quad \partial \Omega_\delta, 
\end{equation}
where $\tilde{\psi}_b, \tilde{\zeta}_b$ are $C^{2,\alpha}$ given functions. 
Then, we have the following 
\begin{theorem}
\label{QPTHM}
For $\gamma\in (-1, + \infty) \backslash  \{0\}$, $0 < \delta_0 << 1$, there exists a $\delta_1 \in (0, \delta_0)$ such that, for any $\delta \in [0, \delta_1)$, 
the Dirichlet problem \eqref{SysPertNew}, \eqref{Dirichlet_Cond_Quasi}
has a solution $$(\psi, \zeta)  \in \big( C^3(\Omega_\delta) \cap C^{2,\alpha}(\bar{\Omega}_\delta)\big)^2.$$
Moreover, the pair 
$(\psi, \zeta)$ is in a neighborwood of $(\vp, 0)$, where $\vp$ is the solution of the Dirichlet problem \eqref{qigual0}
given by Theorem \ref{theorem4.3}. 
\end{theorem}
\begin{proof} 
1. The result follows by applying the Implicit Function Theorem in Banach spaces.
First, we consider the transport equation $(\ref{SysPertNew})_2$.
Indeed, let $\theta \in [0,\delta_0]$ be fixed, and we
assume that the drift $\mathbf{b}= \nabla \psi$ is given,
(where $\psi \in Y$, see below item 2).  
Also we denote for convenience $\omega= \Delta \zeta$. 
Hence we have 
$$
  \mathbf{b} \cdot \nabla \omega= - \omega \, (1 + \dive \mathbf{b}). 
$$
Then, we proceed with the classical method of characteristics, and define the curve $\mathbf{\xi}(r)$ by 
$$
  \left \{
  \begin{aligned}
  \frac{d \xi(r)}{dr}&= \mathbf{b}(\psi(r)), 
  \\[5pt]
  \xi(0)&= \xi_b \quad \text{on \,\, $\big\{ \xi_b \in \partial \Omega_\theta$, such that,  $\big(\mathbf{b} \cdot \nu_\theta\big)(\xi_b)< 0 \big\}$. }  
  \end{aligned}
  \right .
$$
Therefore, denoting $\omega(r)= \omega(\xi(r))$, it follows that 
$$
  \omega(r)= \omega(\xi_b) \exp\big(-\int_0^r (1 + \dive \mathbf{b})(s) \, ds \big), 
$$
where $\omega(\xi_b)= \theta \Delta \tilde{\zeta}_b$. 
Consequently, we have $\Delta \zeta$ in terms of $\psi$ and $\theta \in [0,\delta_0]$. 
Due to a standard argument we obtain also $\nabla^{\perp} \zeta$. 

\medskip
2. Now, we denote  
$$
\begin{aligned}
   &Y_\theta= \big\{ \psi \in C^3(\Omega_\theta) \cap C^{2,\alpha}(\bar{\Omega}_\theta) \,/\,  \psi|_{\partial \Omega_\theta}= \tilde{\psi}_b \big\},   \quad Z_\theta= C^1(\Omega_\theta), 
\\[5pt]   
   &\quad \text{and set} \quad Y= \bigcup_{\theta \in [0,\delta_0]} Y_\theta, \quad Z= \bigcup_{\theta \in [0,\delta_0]} Z_\theta. 
\end{aligned}
$$  
Then, considering equation $(\ref{SysPertNew})_1$
we define the map 
$\mathcal{F}: [0,\delta_0] \times Y \to Z$, given by
$$
\mathcal{F(\theta, \psi)}= c_0^2\,  \Delta \psi - (D^2 \psi) \nabla \psi \cdot \nabla \psi - | \nabla \psi|^2 + 2 \, c_0^2 
- \theta \big( (2 +  \Delta \psi) \, Q_1 + N_1 \big). 
$$
With $(0, \psi_0) \in [0,\delta_0] \times Y$, if 
$$
\begin{aligned}
&(i) \quad \mathcal{F}(0, \psi_0)= 0, 
\\[5pt]
&(ii) \quad D_\psi \mathcal{F}|_{(0,\psi_0)} \quad \text{is a linear homeomorphism of $Y$ in $Z$,} 
\end{aligned}
$$
then we may apply the Implict Function Theorem, and there exists $\delta_1 \in (0,\delta_0)$ and 
a unique mapping $\psi: (0,\delta_1) \to Y$, such that, 
for each $\delta \in (0,\delta_1)$, we have 
$\mathcal{F}(\delta, \psi)= 0$, where $\psi \equiv \psi(\delta) \in Y$ and solve 
$(\ref{SysPertNew})_1$ with $\psi= \tilde{\psi}_b$ on $\partial \Omega_\delta$. 

\medskip
Clearly, applying Theorem \ref{theorem4.3} with $ \psi_0= \vp$ item $(i)$ holds true. 
Therefore, it remains to verify item $(ii)$, which we address below.
We remark that, since $D_\psi \mathcal{F}$
is evaluated at $(0, \psi_0)$, it not necessary to compute the derivative of the term 
$\big( (2 +  \Delta \psi) \, Q_1 + N_1 \big)$. 

\medskip
3. Given $v \in Y$, let us compute 
\begin{eqnarray*}
D_\psi\mathcal{F}\mid_{(0, \psi_0)} [v] 
&=&  \lim_{\tau \to 0} \dfrac{\mathcal{F}(\theta, \psi_0 + \tau v)-\mathcal{F}(\theta, \psi_0)}{\tau}\mathrel{\bigg|}_{(\theta= 0, \psi= \psi_0)}\\
&=& I_1 + I_2 + I_3+I_4+I_5,
\end{eqnarray*}
where at $(0, \psi_0)$, 
\begin{eqnarray*}
I_1&=& \lim_{\tau \to 0} \dfrac{1}{\tau} [c_0^2(\psi_0 + \tau v) \Delta (\psi_0 + \tau v) - c_0^2(\psi_0) \Delta \psi_0],
\\
I_2&=& \lim_{\tau \to 0} \dfrac{1}{\tau} \left[ -(D^2(\psi_0 + \tau v))\nabla (\psi_0 +\tau v) \cdot \nabla (\psi_0 + \tau v) + (D^2 \psi_0) \nabla \psi_0 \cdot \nabla \psi_0 \right],
\\
I_3&=& \lim_{\tau \to 0} \dfrac{1}{\tau} \left[ - |\nabla (\psi_0 + \tau v)|^2 + |\nabla \psi_0 |^2 \right],
\\
I_4&=& \lim_{\tau \to 0} \dfrac{1}{\tau} \left[ -2c_0^2(\psi_0 + \tau v) + 2c_0^2(\psi_0)\right],
\\
I_5&=&  0.
\end{eqnarray*}
To follow, we remember that, 
$$c_0^2 (\psi) = (\gamma-1) \big(-\psi - \dfrac{1}{2}\left| \nabla \psi \right|^2 \big).$$
Therefore, we have 
$$
\begin{aligned}
I_1&=   \lim_{\tau \to 0} \dfrac{1}{\tau} \left[ c_0^2(\psi_0 + \tau v) \Delta (\psi_0 + \tau v) - c_0^2(\psi_0) \Delta \psi_0 \right] \hspace{200pt}
\\[5pt]
&=  \lim_{\tau \to 0} \dfrac{(\gamma -1)}{\tau} \Big( \big( - \psi_0 - \tau v - \frac{1}{2}|\nabla \psi_0|^2 - \tau \nabla \psi_0 \cdot \nabla v - \frac{1}{2} \tau^2| \nabla v|^2 \big) 
\\[5pt]
&\quad \times \quad  (\Delta \psi_0 + \tau  \Delta v) - (- \psi_0 - \frac{1}{2} | \nabla \psi_0|^2 ) \Delta \psi_0 \Big)
\\[5pt]
&=  \lim_{\tau \to 0} \dfrac{(\gamma -1)}{\tau} \left[ (- \psi_0 - \frac{1}{2}| \nabla \psi_0|^2) \tau \Delta v - \tau (v+ \nabla \psi_0 \cdot \nabla v) \Delta \psi_0 \right]
\\[5pt]
&= (\gamma-1) \left[ ( - \psi_0 - \frac{1}{2}| \nabla \psi_0|^2) \Delta v - (v+ \nabla \psi_0 \cdot \nabla v )\Delta \psi_0 \right]
\\[5pt]
&= c_0^2 (\psi_0) \, \Delta v - (\gamma-1) \, \Delta \psi_0( v + \nabla \psi_0 \cdot \nabla v); 
\end{aligned}
$$

$$
\begin{aligned}
I_2&=  \lim_{\tau \to 0} \dfrac{1}{\tau} \big( - (D^2 (\psi_0 + \tau v)) \nabla (\psi_0 + \tau v) \cdot \nabla (\psi_0 + \tau v) + (D^2 \psi_0) \nabla \psi_0 \cdot \nabla \psi_0 \big) \hspace{100pt}
\\[5pt]
&=  \lim_{\tau \to 0} \dfrac{1}{\tau} \Big( - \tau (D^2 \psi_0) \nabla \psi_0 \cdot \nabla v - \tau (D^2 \psi_0) \nabla v \cdot \nabla \psi_0 - \tau^2 (D^2 \psi_0) \nabla v \cdot \nabla v  
\\[5pt]
& \quad - \tau (D^2 v) \nabla \psi_0 \cdot \nabla \psi_0 - \tau^2 (D^2v) \nabla \psi_0 \cdot \nabla v - \tau^2 (D^2v) \nabla v \cdot \nabla \psi_0 - \tau^3 (D^2v) \nabla v \cdot \nabla v \Big)
\\[5pt]
&= - (D^2v) \nabla \psi_0 \cdot \nabla \psi_0 -2 (D^2 \psi_0) \nabla \psi_0 \cdot \nabla v; 
\end{aligned}
$$

$$
\begin{aligned}
I_3&=  \lim_{\tau \to 0} \dfrac{1}{\tau} \left[ - \mid \nabla (\psi_0 + \tau v) \mid^2 + \mid \nabla \psi_0 \mid^2\right] \hspace{300pt}
\\[5pt]
%&=  \lim_{\tau \to 0} \dfrac{1}{\tau} \left[- (\nabla \psi_0 + \tau \nabla v) \cdot ( \nabla \psi_0 + \tau \nabla v) +  \nabla \psi_0 \cdot  \nabla \psi_0 \right]
%\\[5pt]
&=  \lim_{\tau \to 0} \dfrac{1}{\tau} \left[ - 2 \tau (\nabla \psi_0 \cdot \nabla v )- \tau^2 (\nabla v \cdot \nabla v) \right]
\\[5pt]
&= - 2 (\nabla \psi_0 \cdot \nabla v); 
\end{aligned}
$$

$$
\begin{aligned}
I_4&=  \lim_{\tau \to 0} \dfrac{1}{\tau} \left[ -2c_0^2 (\psi_0 + \tau v) + 2c_0^2 (\psi_0) \right] \hspace{300pt}
\\[5pt]
&=  \lim_{\tau \to 0} \dfrac{2(\gamma -1)}{\tau} \left[ \psi_0 + \tau v + \frac{1}{2}(\nabla \psi_0+ \tau \nabla v) \cdot( \nabla \psi_0 + \tau \nabla v) - \psi_0 - \frac{1}{2}(\nabla \psi_0 \cdot \nabla \psi_0 ) \right]
\\[5pt]
&=  \lim_{\tau \to 0} \dfrac{2(\gamma -1)}{\tau} \left[\tau v + \tau (\nabla \psi_0 \cdot \nabla v) + \frac{1}{2} \tau^2 (\nabla v \cdot \nabla v) \right]
\\[5pt]
&= 2 (\gamma -1) \left( v + \nabla \psi_0 \cdot \nabla v \right).
\end{aligned}
$$
Consequently, we obtain 
\begin{equation}
\label{Lv} 
\begin{aligned}
D_\psi\mathcal{F}\mid_{(0, \psi_0)} [v] 
%&= c_0^2(\psi_0) \Delta v - (\gamma-1) \Delta \psi_0 (v+ \nabla \psi_0 \cdot \nabla v) -(D^2v) \nabla \psi_0 \cdot \nabla \psi_0 - 2(D^2 \psi_0) \nabla \psi_0 \cdot \nabla v
%\\
%& \quad  - 2(\nabla \psi_0 \cdot \nabla v) + 2(\gamma-1) (v+ \nabla \psi_0 \cdot \nabla v)\\
&= c_0^2 (\psi_0) \Delta v - (D^2 v) \nabla \psi_0 \cdot \nabla \psi_0 
\\[5pt]
&  - 2(D^2\psi_0) \nabla \psi_0 \cdot \nabla v + [(\gamma-1) (2- \Delta \psi_0) -2] \nabla \psi_0 \cdot \nabla v 
\\[5pt]
& + (\gamma-1)(2-\Delta \psi_0) v
\\[5pt]
& =: \mathcal{L}^{\!\psi_0}[v]. 
\end{aligned}
\end{equation}
That is, $\mathcal{L}^{\!\psi_0}$ is a linear second order partial differential equation in $\Omega_0$, 
which is well defined from the regularity of $\psi_0$.   
Thus we may apply the results in \cite{Gilbarg}, Ch. 6, to show that $\mathcal{L}^{\!\psi_0}$
satisfies $(ii)$. 

\medskip
4. It is not difficult to show that, $\mathcal{L}^{\! \psi_0}$ is continuous. 
To begin, we rewrite 
$$
\begin{aligned}
\mathcal{L}^{\!\psi_0}[v]&= \big(c_0^2 (\psi_0) I_d -  \nabla \psi_0 \otimes \nabla \psi_0 \big) : (D^2 v) 
\\[5pt]
& + \big( (\gamma-1) (2- \Delta \psi_0) I_d -2(D^2\psi_0 + I_d) \big) \nabla \psi_0 \cdot \nabla v 
%&  - 2(D^2\psi_0) \nabla \psi_0 \cdot \nabla v + [(\gamma-1) (2- \Delta \psi_0) -2] \nabla \psi_0 \cdot \nabla v 
\\[5pt]
& + (\gamma-1)(2-\Delta \psi_0) \, v 
\\[5pt]
&= a_{ij}(\xi)  \, \partial^2_{ij} v + b_i(\xi) \, \partial_i v + c(\xi) \, v
\end{aligned}
$$
with obvious notation. 
Let $v_n \to v \in Y$, then $\partial_{ij}^2 v_n \to \partial_{ij}^2 v$ and $\partial_i v_n \to \partial_i v$. Therefore, we have
\begin{eqnarray*}
\lim_{n \to \infty} \mathcal{L}^{\psi_0}[v_n] &=& \lim_{n \to \infty} \big(a_{ij}(\xi) \partial_{ij}^2v_n+ b_i(\xi) \partial_i v_n + c(\xi) v_n \big)
\\
&=& a_{ij}(\xi) \partial_{ij}^2v+ b_i(\xi) \partial_i v + c(\xi) v
\\
&=& \mathcal{L}^{\psi_0} [v].
\end{eqnarray*}

\medskip
5. Here, we prove that $\mathcal{L}^{\psi_0}$ is a bijective operator. 
One remarks that, since there is no sign control on $c(\xi)$, this is a non-trivial equation to be analized, see p. 107 in \cite{Gilbarg}. 
Indeed, we are goind to apply the Fredholm alternative, which could be seen at Theorem 6.15 in that reference. 
It follows from the Fredholm aternative that to prove that $\mathcal{L}^{\psi_0}$ 
is bijective it is sufficient to prove that $\mathcal{L}^{\psi_0}$ is strictly elliptic 
and injective, that is, $\mathcal{L}^{\psi_0}[v]=0$ with $\tilde{\psi}_b \equiv 0$, if and only if, $v=0$.

\medskip
\underline {Claim$\# 1$:} 
 We claim that $\mathcal{L}^{\!\psi_0}$ is uniformly (or strictly) elliptic in $\Omega_0$, that is, there exists $\lambda> 0$, such that
$$
    a_{ij}(\xi) \eta_i \eta_j \ge \lambda |\eta|^2, \quad \text{for all $\xi \in \Omega_0$, and $\eta \in \mathbb{R}^n \setminus \{0\}$.}
$$
 \medskip
Proof of Claim $\# 1$: 1. First, let us show that $\mathcal{L}^{\!\psi_0}$ is elliptic for each $\xi \in \Omega_0$. Indeed, 
it follows that 
$$
\begin{aligned}
a_{ij}(\xi) \eta_i \eta_j &= [c_0^2(\psi_0) - (\partial_1 \psi_0)^2] \eta_1^2 - 2 (\partial_1 \psi_0)(\partial_2 \psi_0) \eta_1 \eta_2 +[c_0^2(\psi_0) - (\partial_2 \psi_0)^2] \eta_2^2
%\\[5pt]
%&=c_0^2(\psi_0) (|\eta|^2-((\partial_1 \psi_0)^2 \xi_1^2 + 2(\partial_1 \psi_0)(\partial_2 \psi_0) \xi_1 \xi_2 + (\partial_2 \psi_0)^2 \xi_2^2)
\\[5pt]
&= c_0^2(\psi_0) \, |\eta|^2 - [(\partial_1 \psi_0) \eta_1+ (\partial_2 \psi_0)\eta_2]^2
\\[5pt]
&= c_0^2(\psi_0) \, |\eta|^2 -  (\nabla \psi_0 \cdot \eta)^2
\\[5pt] 
& \geq (c_0^2(\psi_0) - |\nabla \psi_0|^2) |\eta|^2,
\end{aligned}
$$
where we have used the Cauchy-Schwarz's inequality. 
Since $|\nabla \psi_0(\xi)|^2 < c_0^2(\psi_0(\xi))$ for each $\xi \in \Omega_0$, we obtain 
$$
   a_{ij}(\xi) \eta_i \eta_j >  0
$$
and thus $\mathcal{L}^{\! \psi_0}$ is elliptic.

\medskip
2. At first glance, $\mathcal{L}^{\!\psi_0}$
appears not to be uniformly elliptic because of the sonic line as an intrinsic feature.
Although, we show that applying the Ellipticity Principle Theorem \ref{Ellipticity_principle_theorem}, item $(ii)$: 
there exists $C_0> 0$, such that, for all $\delta> 0$
\begin{equation}
L^2 \le C_0 \delta \quad \mbox{in} \quad \Omega_0.
\end{equation}
Fix $\delta$ such that $C_0\delta <1$. Therefore,
\begin{equation}
L^2 = \dfrac{ | \nabla \psi_0 |^2}{c_0^2} \le C_0 \, \delta,
\end{equation}
and we have 
$$
\begin{aligned}
a_{ij}(\xi) \eta_i \eta_j & \geq c_0^2(\psi_0(\xi)) \; \big(1 - C_0 \; \delta \big) |\eta|^2 
\\[5pt]
&\geq \min_{\xi \in \bar{\Omega}_0} c_0^2(\psi_0(\xi)) \; \big(1 - C_0 \; \delta \big) |\eta|^2,
\end{aligned}
$$
which finish the claim with 
\begin{equation}
\label{Lambda}
  \lambda:=  \min_{\xi \in \bar{\Omega}_0} c_0^2(\psi_0(\xi)) \; \big(1 - C_0 \; \delta \big)> 0. 
\end{equation}

\medskip
Then, we follow the strategy develop in the proof of Theorem \ref{theorem4.1} to show that $\mathcal{L}^{\! \psi_0}$
is injective.
Conveniently, let $\mathcal{L}^{\! \psi_0}_-:= - \mathcal{L}^{\! \psi_0}$. 
Then $\mathcal{L}^{\! \psi_0}$ is injective, if and only if, $\mathcal{L}^{\! \psi_0}_-$ 
is injective, since both operators share the same kernel. Moreover, we write the operator 
$\mathcal{L}^{\! \psi_0}_-$  in divergence form, that is to say, 
$$
   \begin{aligned}
   \mathcal{L}^{\! \psi_0}_-[v]&= - \partial_j(a_{ij}(\xi)  \, \partial_{i} v) +(\partial_j a_{ij}(\xi) - b_i(\xi)) \, \partial_i v - c(\xi) \, v
   \\[5pt]
   &=  - \partial_j(a_{ij}(\xi)  \, \partial_{i} v) + \tilde{b}_i(\xi) \, \partial_i v + \tilde{c}(\xi) \, v, 
   \end{aligned}
$$
where 
$$
  \begin{aligned}
  \tilde{b}_i(\xi)&= \partial_j \big(c_0^2 (\psi_0) \delta_{ij} -  \partial_i \psi_0 \; \partial_j \psi_0 \big)
  \\[5pt]
  &\quad -\big( (\gamma-1) (2- \partial^2_{kk} \psi_0) \delta_{ij} -2(\partial^2_{ij} \psi_0 + \delta_{ij}) \big) \partial_j \psi_0, 
  \\[5pt]
  \tilde{c}(\xi)&= - (\gamma-1)(2-\Delta \psi_0). 
  \end{aligned}
$$
Then, we consider the bilinear form $\tilde{B}[v,w]$ associated with the shifted operator $\mathcal{L}^{\! \psi_0}_- + \kappa_0 I$, 
for some $\kappa_0> 0$ to be chosen a posteriori, defined by 
\begin{equation}
\label{bilinear_form}
    \tilde{B}[v,w]:= \int_{\Omega_0} a_{ij}(\xi) \;  \partial_i v \; \partial_j w \, d\xi 
    + \int_{\Omega_0} \tilde{b}_i(\xi) \, \partial_i v \, w \; d\xi 
    + \int_{\Omega_0} (\tilde{c} + \kappa_0) \; v \; w \; d\xi, 
\end{equation}
which is well-defined for instance for $v, w \in H^1_0(\Omega_0)$. 

\medskip
\underline {Claim$\# 2$:} 
 We claim that the bilinear form $\tilde{B}[\cdot,\cdot]$ is 
 continuous, bounded and coercive in $H^1_0(\Omega_0)$. 
 Moreover, the linear operator 
 $\mathcal{L}^{\! \psi_0}_-$ is injetive. 
 
  \medskip
Proof of Claim $\# 2$: 
The proof is analogous to the arguments used in items 2 and 3 of the 
proof of Theorem \ref{theorem4.1}, with no essential modifications, 
and hence it is omitted.

\medskip
6. Finally, since $\mathcal{L}^{\psi_0}$ is a continuous linear and bijective operator, 
it follows from the Open Mapping Theorem that $\mathcal{L}^{\psi_0}= D_\psi \mathcal{F}|_{(0,\psi_0)}$ is 
a linear homeomorphism of $Y$ in $Z$, see Corollary 2.7 in Brezis \cite{BREZIS}. 
Consequently, the proof of Theorem \ref{QPTHM} is complete.
\end{proof}

%%%%%%%%%%%%%%%%%%%%
\section*{Data availability statement}
%%%%%%%%%%%%%%%%%%%%
Data sharing is not applicable to this article as no data sets were 
generated or analysed during the current study. 

%%%%%%%%%%%%%%%%%%%%%%%%%%%%%%
 \section*{Acknowledgements}
%%%%%%%%%%%%%%%%%%%%%%%%

The author Wladimir Neves has received research grants from CNPq
through the grant  313005/2023-0, 403675/2025-1,
and also by FAPERJ 
(Cientista do Nosso Estado) through the grant E-26/204.171/2024.

\end{document}